\documentclass[preprint,12pt]{elsarticle}

\usepackage{algorithm}
\usepackage{algorithmic}

\usepackage{mathtools}
\usepackage{bm} 
\usepackage{dutchcal}
\usepackage{graphicx}               
\usepackage{subfigure}
\usepackage{caption}
\usepackage[section]{placeins}
\usepackage{float}

\usepackage{amssymb}

\usepackage{lineno}

\biboptions{numbers,sort&compress}

\date{}

\journal{}

\begin{document}

\begin{frontmatter}




\title{ 
A stabilized mixed finite element scheme for frictional contact mechanics and shear failure analyses in deformable media with crossing fractures
\tnoteref{t1,t2}}





\author[tudewi]{Luyu Wang \corref{cor1}} 
\ead{wang.luyu@cnrs.fr; luyu.wang@hotmail.com}

\author[tudewi]{Cornelis Vuik }
\ead{c.vuik@tudelft.nl}

\author[tudgitc]{\mbox{Hadi Hajibeygi}}
\ead{h.hajibeygi@tudelft.nl}

\address[tudewi]{Department of Applied Mathematics, Delft University of Technology, 2628 CD Delft, the Netherlands}
\address[tudgitc]{Department of Geoscience and Engineering, Delft University of Technology, 2628 CV, Delft, the Netherlands}

\cortext[cor1]{Corresponding author. Current affiliation is GeoRessources Lab., CNRS, 54500 Nancy, France. }


\begin{abstract}

Simulation of frictional contact and shear failure of fractures in fractured media is of paramount important in computational mechanics. In this work, a preconditioned mixed-finite element (FE) scheme with Lagrange multipliers is proposed in the framework of constrained variational principle, which has the capability to handle frictional contact and slip of the multiple crossing fractures. 
The slippage, opening and contact traction on fractures are calculated by the resulted saddle-point algebraic system. 
A novel treatment is devised to guarantee physical solutions at the intersected position of crossing fractures. A preconditioning technique is introduced to re-scale the resulting saddle-point algebraic system, to preserve the robustness of the system. 
An iteration strategy, namely monolithic-updated contact algorithm, is then designed to update the two primary unknowns (displacement and Lagrange multiplier) in one algebraic block. 
Then, a series of numerical tests is conducted to study the frictional contact and shear failure of single- and multi-crossing fractures. Benchmark study is performed to verify the presented mixed-FE scheme. Two tests with crossing fractures are studied, in which the slippage and opening can be calculated. The effects of crossing fractures on the deformation field are observed in the simulation, in which the variation of slippage, opening and stress intensity factor are analyzed under different loading conditions.

\end{abstract}


%
%



\end{frontmatter}


\section{Introduction}
\label{section:SecIntroduction}

The subject of computational contact mechanics is a pillar in the field of applied and computational mechanics \cite{Barber2000,Lorenzis2017}. Especially, it plays an essential role in geoscience applications, such as reservoir engineering, subsurface energy exploration, nuclear waste disposal and induced seismicity (i.e., faults activation)  \cite{Jha2014, Garipov2016,Aagaard2013}. Discrete fractures widely exist in geological fields with significantly contrasting hydraulic and mechanical properties. The formation of fracture networks and their evolution depend on the manner of loading condition imposed by the surrounding geological environment \cite{Adler2013,Berkowitz2002,Wang2019b,Wang2020}. Moreover, the multiple crossing fractures lead to several complicated mechanical behaviors, specially when the fractures are subjected to compressive loading. In this situation, the two sides of fracture plate contact and interact with each other, which results in contact tractions. This is totally different from the scenario of tensile loading, as in the situation of hydraulic-driven fractures, in which the contact constraints disappear as the two sides of the fractures are separate \cite{Secchi2012,Lecampion2018}. Mathematically, fractures are viewed as "internal boundaries" inside the interested fractured domain, with respect to the "external boundary" being the outline of the entire domain \cite{Sofonea2012,Wriggers2006,Yastrebov2013,Wang2019a}. However, several challenges appear in the complex situation of crossing fractures \cite{Annavarapu2013}, in which proper treatments at the intersected position of crossing fractures need more attentions. 

The classical computational contact mechanics handles the contact and impact phenomenon between two separate deformable (or rigid) bodies \cite{Wriggers2006,Zhong1993,Meguid2008}. In contrast to this, a situation that the fractures embedded in a fractured medium frequently occurs in geoscience applications \cite{Franceschini2016,Zhao2021,Ren2021}. Moreover, the discrete fractures induce the multi-contact surfaces, which is an unusual configuration in classical contact mechanics, but quite common in geoscience applications. 
To this end, the classical theory is extended to the category of fault mechanics \cite{Garipov2016,Aagaard2013,Jha2014,Franceschini2016,Zhao2021}, which requires the capability to calculate slippage and opening of the frictional contact behavior. On the other hand, to simplify the computational workflow, sometimes the explicit treatment of contact mechanics is neglected \cite{Ren2021}. This simplification is feasible in case of hydro-fracturing problem, when fractures are under tensile forces \cite{Secchi2012,Lecampion2018}.
For compressive force, however, such simplifications lead to significant errors. The particular interest is to develop a computational model which allows for explicit consideration of both opening and slippage (i.e. normal and tangential displacements). In this work, we include frictional contact mechanics and sliding by Lagrange multipliers, in which both the slippage and opening can be accurately modeled. 

In the framework of continuum mechanics, the category of mesh-like numerical methods is applied and improved to model the frictional contact and sliding. Among these methods, finite element method (FEM) is the most widely used method in computational contact mechanics \cite{Lorenzis2017,Jha2014,Garipov2016,Aagaard2013,Wriggers2006,Yastrebov2013,Zhong1993,Meguid2008}. Recently, the extended finite element method (XFEM) \cite{Liu2010,Liu2008} and finite volume method (FVM) \cite{Berge2020,Ucar2018} have been proposed to handle contact in fractured media. However, the classical contact mechanics studies the single contact surface. Based on this, the numerical methods for multi-contact surfaces in fractured media have been devised. In the aspect of treatment of contact constraints, penalty method \cite{Mergheim2004,Ferronato2008} and Lagrange multiplier method \cite{Pietrzak1999,Simo1992} are the two main schemes originated from the variational principle, which renders the fundamental of FEM. 
Meanwhile, an FVM scheme, which considers Lagrange multiplier as an additional unknown, combined with the variationally consistent hybrid discretization has been proposed in the literature \cite{Berge2020}, in which the multi-point stress approximation was used to discretize the governing equations. XFEM has been developed by the penalty method \cite{Liu2008} and Lagrange multiplier method \cite{Liu2010} for modeling contact. 
The key idea of penalty method is to apply the springs with large stiffness to connect the two sides of contact surfaces, thus one does not introduce the additional unknowns \cite{Mergheim2004,Ferronato2008}. The drawback is the resulted linear algebraic system is ill-conditioned. Furthermore, the crucial point is that the penalty method is inaccurate since it enforces springs with large stiffness values. The Lagrange multiplier method, however, introduces an additional unknown \cite{Simo1992,Simo1985}, namely the Lagrange multiplier. It has an underlying physical meaning, i.e. the traction on a contact surface. Its inconvenience is the additional cost of solving the augmented algebraic system with the so-called saddle-point structure \cite{Franceschini2019,Auricchio2017,Pestana2015}. Moreover, the additional unknowns expand the size and disturb the sparsity of the resulted linear systems, so the development of efficient linear solvers become even more important. More precisely, the Lagrange multiplier method results in the mixed-FE scheme \cite{Cescotto1993,Franceschini2019,Franceschini2019}, in which the displacement and Lagrange multiplier can be calculated at one single global algebraic block. Other methods, such as the Nitsche's method \cite{Hansbo2005,Annavarapu2012}, augmented Lagrange multiplier method \cite{Simo1992,Nejati2016} and mortar method \cite{Farah2015,Seitz2016} have been also proposed. In fact, all of these methods are based on the basic principle of penalty method and Lagrange multiplier method in the framework of variational principle. 
The high nonlinear property of the contact system is captured by the so-called Karush-Kuhn-Tucker (KKT) condition \cite{Barber2000,Lorenzis2017,Sofonea2012,Wriggers2006,Yastrebov2013}, which is integrated into the standard variational principle and resulting in the constrained variational principle. The prospective directions on computational contact mechanics is focused on the treatment of multi-fractured systems, including crossing fractures, and development of stable \cite{Pestana2015,Franceschini2019,VorstVuik2017} and scalable iterative solvers \cite{Franceschini2019} for the resulting saddle-point systems. Such developments would allow for modeling large-scale systems within the industrial applications. 

Despite being crucially important, accurate and stable modeling of crossing multi-fractured media remains to a large extend unexplored. This work develops a novel approach to resolve this limitation. First the FE-based system on an unstructured mesh is developed in which fractures are confined at the matrix element interfaces. Then a constraint is introduced to guarantee physical solutions in presence of crossing fractures. The resulting system contains displacement unknowns and Lagrange multipliers, leading to the so-called saddle-point. We develop a novel scaling algorithm to improve the system condition number and thus leads to a robust solution strategy. The developed method is being benchmarked against analytical methods, and tested for several proof-of-concept numerical test cases.


This paper is structured as follows. First, the formulation of frictional contact and shear failure on multiple crossing fractures is presented in Section \ref{section:SecMathFormulation}. Then, the contact constraints are integrated into variational principle using Lagrange multiplier. Galerkin FE approximation is applied in Section \ref{section:SecDiscretization} to discretize the system. In Section \ref{section:SecSolutionStrategy}, the solution strategy is proposed with the devised monolithic-updated contact algorithm, in which the precondition is used to handle the resulted saddle-point system. Finally, a series of numerical tests is performed in Section \ref{section:NumericalTests} to verify the proposed method and to show the contact behavior on multiple fractures.


\section{Physical model}
\label{section:SecMathFormulation}

The formulation is presented after a revisit on the theory of contact mechanics. Both the isolated fracture and crossing fractures can be handled in this framework. The compressive state would lead to the activation of contact constraints, while the tensile state remains the standard scheme of elasticity.


\subsection{Contact model of multi-fractures}
\label{subsection:ContactModel}

\begin{figure}[H]
\centering
\includegraphics[width=9cm]{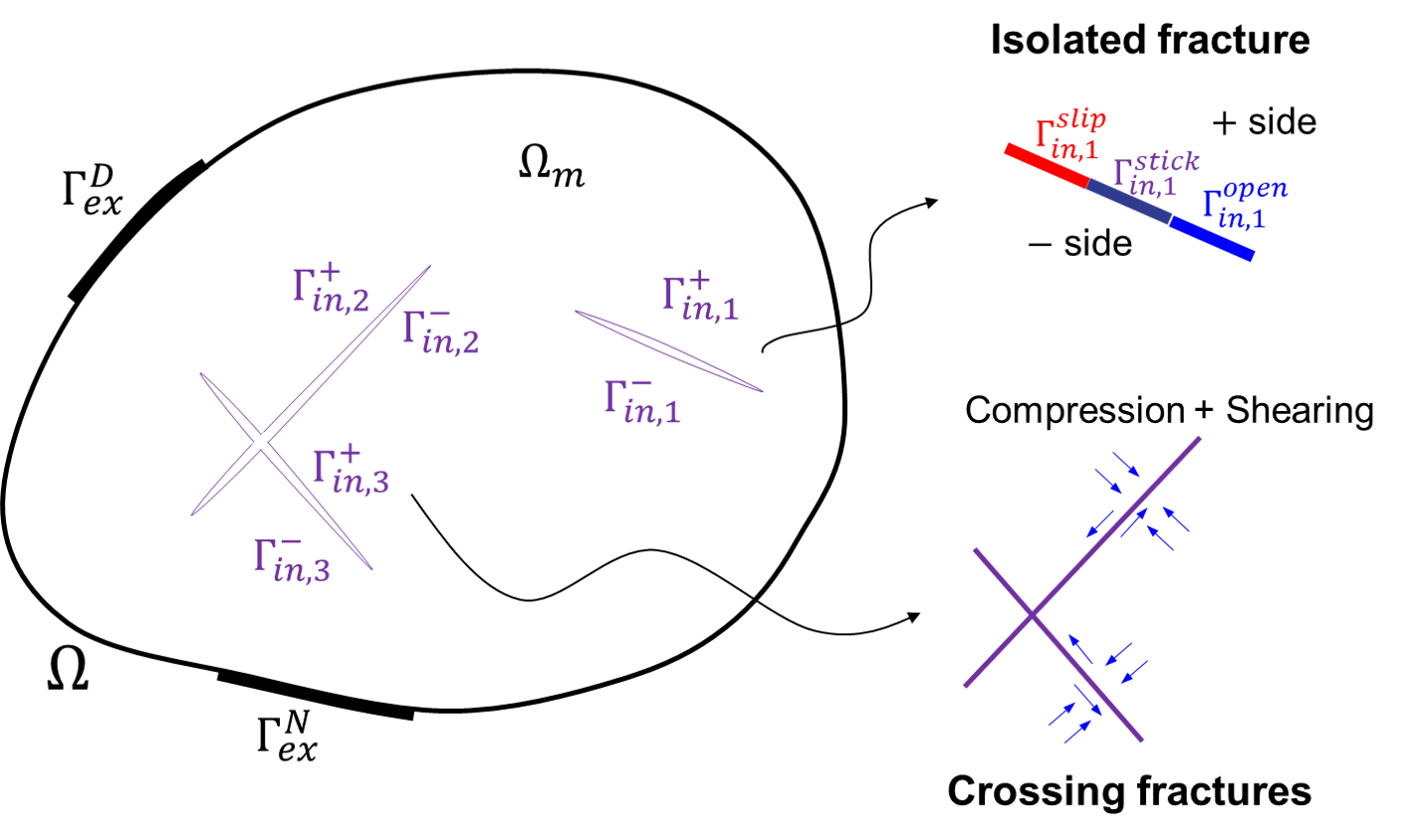}
\caption{ Contact model of multi-fractures containing crossing fractures. Each fracture has a positive and a negative surfaces. Each contact surface can be divided into different portions with state of stick/slip/open. }
\label{fig:fig_ContactModel}
\end{figure}

The porous medium is denoted as $\Omega^m$. 
As shown in Fig. \ref{fig:fig_ContactModel}, a set of the discrete fractures $\omega$ is distributed in the domain and is modeled as the internal boundaries $\Gamma_{in}$ inside $\Omega^m$. 
Each of the fracture surfaces is decomposed into the positive side $\Gamma_{in}^{+}$ and the negative side $\Gamma_{in}^{-}$, such that $\Gamma_{in} = \Gamma_{in}^{+} \cup \Gamma_{in}^{-}$. 
For multiple fractures, it reads: 
\begin{equation}
    \omega := \Gamma_{in} = \sum_{i=1}^{N^f} \Gamma_{in,i} =  \sum_{i=1}^{N^f} \left( \Gamma_{in,i}^{+} \cup \Gamma_{in,i}^{-} \right)
\end{equation}
where $N^f$ is the number of fractures. According to the convention in computational contact mechanics \cite{Lorenzis2017,Sofonea2012,Wriggers2006,Yastrebov2013}, we also adopt the terms of master surface $\Gamma_{in}^{-}$ and slave surface $\Gamma_{in}^{+}$.

The external boundary $\Gamma_{ex}$ is expressed by $\Gamma_{ex} = \Gamma_{ex}^{D} \cup \Gamma_{ex}^{N} \;\; {\rm and} \;\; \Gamma_{ex}^{D} \cap \Gamma_{ex}^{N} = \emptyset$, with Dirichlet- and Neumann- boundaries, namely $\Gamma_{ex}^{D}$ and $\Gamma_{ex}^{N}$. 
The domain $\Omega$ can be decomposed into $\Omega = \Omega^m \cup \Gamma_{in}$, with the boundaries of external- and internal- types $\partial \Omega = \Gamma_{ex} \cup \Gamma_{in}$.

One of the novelties in this work is to model the contact behavior and sliding of $\Gamma_{in}$. The constraint on contact surfaces depends on the type of loading as well as the frictional law, which controls the slip on fractures. Especially, the Karush-Kuhn-Tucker (KKT) condition \cite{Barber2000,Lorenzis2017,Sofonea2012,Wriggers2006} would be active once the contact surface is imposed by compression. It is used to confine the unphysical effects of surface penetration and mesh overlapping, as displayed in Fig. \ref{fig:fig_ContactStates}. The details of contact constraint will be introduced in Section \ref{subsection:Fractures}.


\subsection{Governing equation of the elastic medium}
\label{subsection:ElasticMedium}

The host matrix is assumed as an elastic medium. The assumptions of infinitesimal deformation and the quasi-static contact are employed. The deformation of $\Omega$ is captured by the momentum balance combining with Hooke's law: 
\begin{equation}
\nabla \cdot \left( \mathbb{C} : \nabla^s \bm{u} \right)
+ \bm{f} = \bm{0}
\quad    \rm{on}  \; \Omega
\label{eq:MomentumBalance}
\end{equation}
with the body force $\bm{f}$ and elasticity tensor $\mathbb{C}$. $\bm{u}$ is the displacement vector, $\nabla^s$ the symmetric gradient.

It is straightforward to define boundary conditions for the external boundary. The pre-defined displacement $\bm{\bar{u}}^{ex}$ and traction $\bm{\bar{t}}^{ex}$ are given as: 
\begin{equation}
\begin{aligned}
\bm{u} = \bm{\bar{u}}^{ex}  
\quad \rm{on} \; &\Gamma_{ex}^{D}  \\
\bm{\sigma} \cdot \bm{n}^{ex} = \bm{\bar{t}}^{ex} 
\quad  \rm{on} \; &\Gamma_{ex}^{N}  \\
\end{aligned}
\label{eq:ExternalBC}
\end{equation}
with the outward unit vector $\bm{n}^{ex}$ which points outward of the external boundary.

The situation is totally different once fractures are introduced. Providing a reasonable and correct condition is complicated for $\Gamma_{in}$ (also known as the fractures) because of the nonlinear constraints on fractures \cite{Lorenzis2017,Sofonea2012}. The system of nonlinear inequalities (the KKT condition) \cite{Lorenzis2017,Wriggers2006} are introduced: 
\begin{equation}
\bm{C} \left( \bm{u}^{f} , \bm{t}^{f} \right) \geqslant \bm{0}
\quad  \rm{on} \; \Gamma_{in} 
\label{eq:NonlInearinequalities}
\end{equation}
with the displacement $\bm{u}^{f}$ and traction $\bm{t}^{f}$ on the contact surface. The expanded form will be introduced in Section \ref{subsection:Fractures}.

\begin{figure}[H]
\centering
\includegraphics[width=9.5cm]{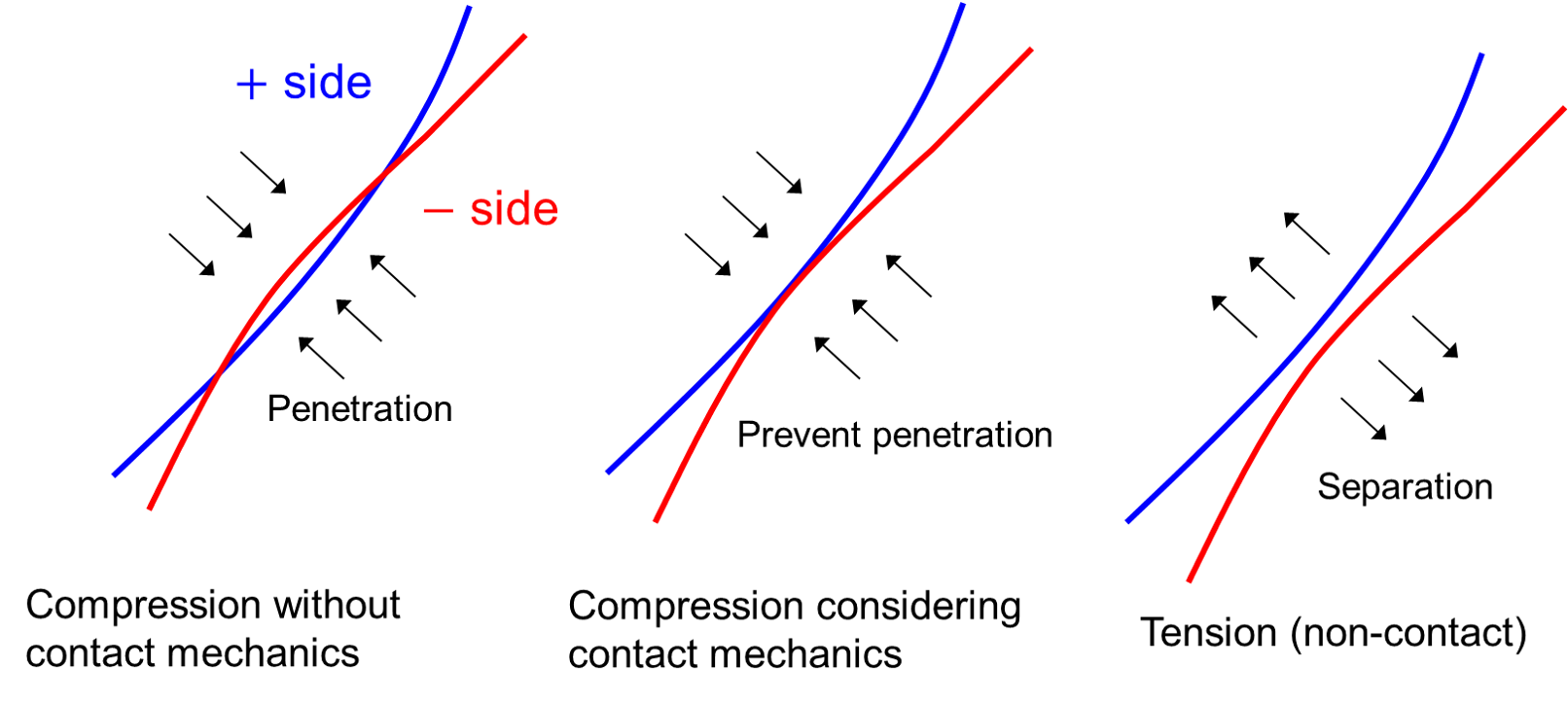}
\caption{ Contact states under different loading. When the contact surface is imposed by tension, the contact constraints can be neglected. If a compression is applied, the constraint condition related to contact mechanics is active. Otherwise, two surfaces would be penetrated. }
\label{fig:fig_ContactStates}
\end{figure}

The external boundary condition Eq. (\ref{eq:ExternalBC}) is allowed to be compressive or tensile. The later case can be treated easily by the standard finite element framework. The fracture surfaces exhibit complex behaviors if compression is applied, as shown in Fig. \ref{fig:fig_ContactStates}.

\subsection{Constraints on local system of fracture}
\label{subsection:Fractures}

The system of additional constraint equations Eq. (\ref{eq:NonlInearinequalities}) on fractures is taken as a complement of the standard boundary condition Eq. (\ref{eq:ExternalBC}). Generally speaking, as shown in Fig. \ref{fig:fig_ContactModel}, a certain fracture can be slip, stick or open. Furthermore, the contact surface is divided into several portions $\Gamma_{in} = \Gamma_{in}^{stick} \cup \Gamma_{in}^{slip} \cup \Gamma_{in}^{open}$.

\begin{figure}[H]
\centering
\includegraphics[width=4.5cm]{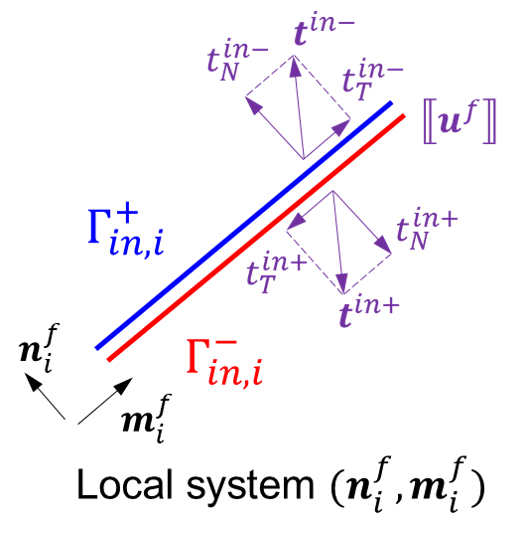}
\caption{ The local system attached to a certain fracture. }
\label{fig:fig_LocalSystem}
\end{figure}

The local coordinate system $\left( \bm{n}^f, \bm{m}^f \right)$, defined by the unit-normal vector $\bm{n}^f$ and the unit-tangential vector $\bm{m}^f$, is attached to each of fractures, as shown in Fig. \ref{fig:fig_LocalSystem}. 
To define the sign of $\bm{n}^f$, we classify the positive and negative vectors, $\bm{n}^{+}$ and $\bm{n}^{-}$, for the two sides $\Gamma_{in}^{+}$ and $\Gamma_{in}^{-}$, respectively. 
The default directions of unit vector and traction used in the presented numerical scheme are defined:
\begin{equation}
\bm{n}^f = \bm{n}^{-} = -\bm{n}^{+} 
\;\; {\rm and} \;\;
\bm{t}_N^{f} = \bm{t}_N^{in-} = -\bm{t}_N^{in+}
\label{eq:VectorTraction}
\end{equation}
where $\bm{n}^f$ is a vector that points from the negative side to the positive side. $\bm{t}_N^{in\pm}$ is the normal component of contact traction on $\pm$ sides. 

The unit vector $\bm{n}^f$ is used to define the normal component of a local system. We adopt the assumption of infinitesimal deformation, thus $\bm{n}^f$ would not be changed at the update time $t+\Delta t$: 
\begin{equation}
 \left. \bm{n}^{f} \right \vert ^{t} =  \left. \bm{n}^{f} \right \vert ^{t+\Delta t} 
\end{equation}

The normal traction of contact $\bm{t}_N^{f}$ is calculated by stress tensor $\bm{\sigma}$ based on the traction vector $\bm{t}^{f}$, which is given by: 
\begin{equation}
\bm{t}^{f} = \bm{\sigma} \cdot \bm{n}^{f} = {t}_N^{f} \bm{n}^f + {t}_T^{f} \bm{m}^f
\label{eq:TractionVector}
\end{equation}
where the sign of compression is negative, vise versa.

Fig. \ref{fig:fig_LocalSystem} depicts a local coordinate system $\left( \bm{n}^f, \bm{m}^f \right)$ for a certain fracture. The quantity $[\![ \bm{u}^f ]\!]$ on a contact surface is introduced to evaluate the relative deformation between the two sides $\Gamma_{in}^{+}$ and $\Gamma_{in}^{-}$: 
\begin{equation}
[\![ \bm{u}^f ]\!] = \bm{u}^{f+} - \bm{u}^{f-}
\label{eq:RelativeDisp}
\end{equation}
with the absolute displacements $\bm{u}^{f+}$ and $\bm{u}^{f-}$ on two sides of the contact surface. The two components of relative displacement are defined as: 
\begin{equation}
[\![ {u}^f_N ]\!] = [\![ \bm{u}^f ]\!] \cdot \bm{n}^f \;\; {\rm and} \;\; 
[\![ {u}^f_T ]\!] = [\![ \bm{u}^f ]\!] \cdot \bm{m}^f
\label{eq:ComponentRelativeDisp}
\end{equation}

In this way, the additional constraint conditions are used to confine the mechanical behavior of contact surface on the normal direction: 
\begin{equation}
\begin{aligned}
g_N = [\![ {u}^f_N ]\!]  \geqslant	0 
\quad  \rm{on} \; &\Gamma_{in}  \\
{t}_N^{f} \leqslant	0
\quad  \rm{on} \; &\Gamma_{in}  \\
g_N {t}_N^{f} = 0
\quad  \rm{on} \; &\Gamma_{in}  \\
\end{aligned}
\label{eq:KKTCondition_Normal}
\end{equation}
and the tangential direction: 
\begin{equation}
\begin{aligned}
\tau_c = c - {t}_N^{f} \tan \varphi 
\quad  \rm{on} \; &\Gamma_{in}  \\
{t}_T^{f} = \tau_c \frac{[\![ {u}^f_T ]\!]}{\Vert [\![ {u}^f_T ]\!] \Vert_2}
\quad  \rm{on} \; &\Gamma_{in}  \\ 
\end{aligned}
\label{eq:KKTCondition_tangential}
\end{equation}
where symbol $\Vert \cdot \Vert_2$ is the 2-norm. The critical value of $\tau_c$ is determined by Mohr-Coulomb criterion. $c$ and $\varphi$ are parameters to control frictional behavior. The sign of ${t}_T^{f}$ is of great important to calculate the dissipate energy induced by contact friction. In this work, we employ the so-called maximum plastic dissipate principle \cite{Wriggers2006,Simo2006,Franceschini2016} to obtain the direction of ${t}_T^{f}$, in which $[\![ {u}^f_T ]\!]$ is the weighted quantity.

Eqs. (\ref{eq:KKTCondition_Normal}) and (\ref{eq:KKTCondition_tangential}) are the expanded forms of Eq. (\ref{eq:NonlInearinequalities}), which hold true on contact surface of each fracture. They provide a set of constraint equations to define the frictional contact, which are considered as requisite complements to Eqs. (\ref{eq:MomentumBalance}) and (\ref{eq:ExternalBC}).

\section{Numerical discretization}
\label{section:SecDiscretization}

The constraints of frictional contact are integrated into the mixed-finite element scheme through Lagrange multipliers. Then, the Galerkin finite element approximation is used to derive the fully discretized scheme resulting an unified computational formula.


\subsection{The variational principle for contact of multiple fractures}

In this section, the variational principle is extended to include the contact mechanics and sliding of multi-fractures. 
Note that Eq. (\ref{eq:KKTCondition_tangential}) is inactive if the current traction component ${t}_T^{f}$ is lower than the critical traction $\tau_c$. 

We use the Galerkin finite element formulation (GFEM) \cite{Wriggers2006,Zienkiewicz2000} to integrate the contact constraints Eqs. (\ref{eq:KKTCondition_Normal}) and (\ref{eq:KKTCondition_tangential}). The standard variational principle (SVP) is extended to a more general case. According to the Lagrange multiplier method \cite{Bertsekas1982,Ito2008,Zienkiewicz2000}, SVP is generalized to the constrained variational principle (CVP). The total energy functional is introduced as: 
\begin{equation}
\Pi^* \left( \bm{u},  \bm{u}^f ,\bm{\lambda} \right) 
= \Pi_u \left(\bm{u} \right)
+ \Pi^{CL} \left( \bm{u}^f ,\bm{\lambda} \right)
\label{eq:Functional}
\end{equation}
with the elastic functional $\Pi_u \left(\bm{u} \right)$ defined by SVP, which is related to the unknown (displacement $\bm{u}$). The novelty is the contact functional $\Pi^{CL} \left( \bm{u}^f ,\bm{\lambda} \right)$ defined by CVP through Lagrange multiplier $\bm{\lambda}$.

The Lagrange multiplier vector reveals an underlying meaning in physical aspect. The component form  $\bm{\lambda} = \left[ \lambda_N \;\; \lambda_T \right]^T$ indicates that it equals the components of contact traction on the surface. 
To solve the system, the first-order variation of Eq. (\ref{eq:Functional}) is written as: 
\begin{equation}
\delta \Pi^* \left( \bm{u},  \bm{u}^f ,\bm{\lambda} \right) 
= \delta \Pi_u \left(\bm{u} \right)
+ \delta \Pi^{CL}_{u} \left( \bm{u}^f ,\bm{\lambda} \right) 
+ \delta \Pi^{CL}_{\lambda} \left( \bm{u}^f ,\bm{\lambda} \right)
\label{eq:Variation1orderFunctional}
\end{equation}

The last two terms in right hand side of Eq. (\ref{eq:Variation1orderFunctional}) are induced by the effects from frictional contact and sliding. 
The Lagrange multiplier method provides a way to integrate contact constraints with the help of Eqs. (\ref{eq:MomentumBalance}), (\ref{eq:ExternalBC}) and (\ref{eq:NonlInearinequalities}): 
\begin{small}
\begin{equation}
\begin{aligned}
\delta \Pi^* \left( \bm{u},  \bm{u}^f ,\bm{\lambda} \right) 
&= 
\int_{\Omega} \delta \bm{\varepsilon}^T \bm{\sigma} d\Omega
- \int_{\Omega} \delta \bm{u}^T \bm{f} d\Omega 
- \int_{{{\Gamma}}_{ex}^N} \delta \bm{u}^T \bm \bar{\bm{t}}^{ex} d\Gamma \\
&+ 
\sum_{i=1}^{N^f} \left(
\int_{\Gamma_{in,i}} \delta \bm{\lambda}^T \bm{C} \left( \bm{u}^{f} , \bm{t}^{f} \right) d\Gamma 
+ \int_{\Gamma_{in,i}}  \bm{\lambda}^T \delta \bm{C} \left( \bm{u}^{f} , \bm{t}^{f} \right) d\Gamma
\right)
\end{aligned}
\label{eq:LagrangeMuliterMethod}
\end{equation}
\end{small}

The variational principle requires $\delta \Pi^* \left( \bm{u},  \bm{u}^f ,\bm{\lambda} \right) = 0$ with respect to displacement $\bm{u}$: 
\begin{small}
\begin{equation}
\begin{aligned}
\int_{\Omega} \delta \bm{\varepsilon}^T \bm{\sigma} d\Omega
- \int_{\Omega} \delta \bm{u}^T \bm{f} d\Omega 
- \int_{{{\Gamma}}_{ex}^N} \delta \bm{u}^T \bm \bar{\bm{t}}^{ex} d\Gamma 
 + 
\sum_{i=1}^{N^f} 
\int_{\Gamma_{in,i}}  \bm{\lambda}^T \delta \bm{C} \left( \bm{u}^{f} , \bm{t}^{f} \right) d\Gamma
= 0
\end{aligned}
\label{eq:VariationalPrinciple_Disp}
\end{equation}
\end{small}
and with respect to Lagrange multiplier $\bm{\lambda}$: 
\begin{small}
\begin{equation}
\begin{aligned}
\sum_{i=1}^{N^f} 
\int_{\Gamma_{in,i}} \delta \bm{\lambda}^T \bm{C} \left( \bm{u}^{f} , \bm{t}^{f} \right) d\Gamma 
  = 0
\end{aligned}
\label{eq:VariationalPrinciple_Lambda}
\end{equation}
\end{small}
with the number of fractures $N^f$. The situation of multi-fractures is then included.

\subsection{Weak forms for stick and slip states}
\label{subsection:SubSecWeakForms}

The first three terms in Eq. (\ref{eq:VariationalPrinciple_Disp}) can be easily treated by the standard Galerkin finite element method. The key point is the fourth term at the left hand side, which describes the effects from contact and sliding. Assuming the initial state of contact surface is stick, the weak form can be expressed. We refer to \ref{appendix_weakForms}.

The situation of slip state is different from stick state once the current tangential component of contact stress reaches the critical stress $\tau_c$, as indicated in Eq. (\ref{eq:KKTCondition_tangential}). In this case, the tangential component is obtained by Mohr-Coulomb criterion, which means $\lambda_T = \tau_c$. We define the indicator function for slip contact: 
\begin{equation}
\begin{aligned}
sign =  \frac{[\![ {u}^f_T ]\!]} {\Vert [\![ {u}^f_T ]\!] \Vert_2} = 
\begin{cases}
1    \quad   \text{ slip along direction of } \bm{m}^f\\
-1   \quad   \text{ slip along direction of } -\bm{m}^f
\end{cases}
\end{aligned}
\label{eq:IndicatorFunction}
\end{equation}

Consequently, the component $\lambda_T$ of unknown $\bm{\lambda}$ has been already known in slip contact, expressed by $\lambda_T = \left( c - \lambda_N \tan \varphi \right) sign$. 
In this way, $\lambda_T$ is connected to the unknown $\lambda_N$. 
Furthermore, compared to the system of stick contact Eq. (\ref{eq:WeakFormState}), the weak form of slip contact needs to be rewritten. We refer to \ref{appendix_weakForms} for the detail.

The slip occurs if Eqs. (\ref{eq:KKTCondition_Normal}) and (\ref{eq:KKTCondition_tangential}) are active (compression). In contrast, if the contact surface is subjected to tension, the situation (open state) is much simpler than compression. 
The weak form would be reduced to the SVP framework, which can be expressed by the first three terms at right hand side of Eq. (\ref{eq:LagrangeMuliterMethod}).

\subsection{Galerkin finite element approximations}
\label{subsection:SubSecGalerkinFEM}

Eqs. (\ref{eq:WeakFormState}), (\ref{eq:WeakFormSlip_Disp}) and (\ref{eq:WeakFormSlip_Lambda}) provide a general formulation which captures the contact and shear failure of multi-fractured media with different types of loading, including compressive, tensile and shear forces. The weak forms can be fully discretized using Galerkin finite element method \cite{Wriggers2006,Zienkiewicz2000}. 

As the requirement of finite element method, the solutions of $\bm{\lambda}$ and $\bm{u}$ would be subjected to the square integrable functional spaces \cite{Brenner2007}. The spaces of displacement and Lagrange multiplier are defined, and the test functions are restricted to Sobolev spaces \cite{Brenner2007,Franceschini2019,Franceschini2020}. $\bm{\lambda}$ and $\bm{u}$ are approximated as: 
\begin{equation}
\begin{aligned}
\bm{u} \approx \bm{u}_h = \sum_{i=1}^{n_{node}} \bm{N}^u_i \bm{U}_i
= \bm{N}^u \bm{U}
\\
\bm{\lambda} \approx \bm{\lambda}_h = \sum_{j=1}^{n_{cp}} \bm{N}^{\lambda}_j \bm{\Lambda}_j
= \bm{N}^\lambda \bm{\Lambda}
\end{aligned}
\label{eq:PrimaryApproximations}
\end{equation}
with displacement vector $\bm{U}$ and Lagrange multiplier vector $\bm{\Lambda}$ at grid vertices, and the shape function matrices $\bm{N}^u$ and $\bm{N}^\lambda$. Note that the numbers of nodes and contact pairs are denoted by $n_{node}$ and $n_{cp}$.

\begin{figure}[H]
\centering
\includegraphics[width=10.5cm]{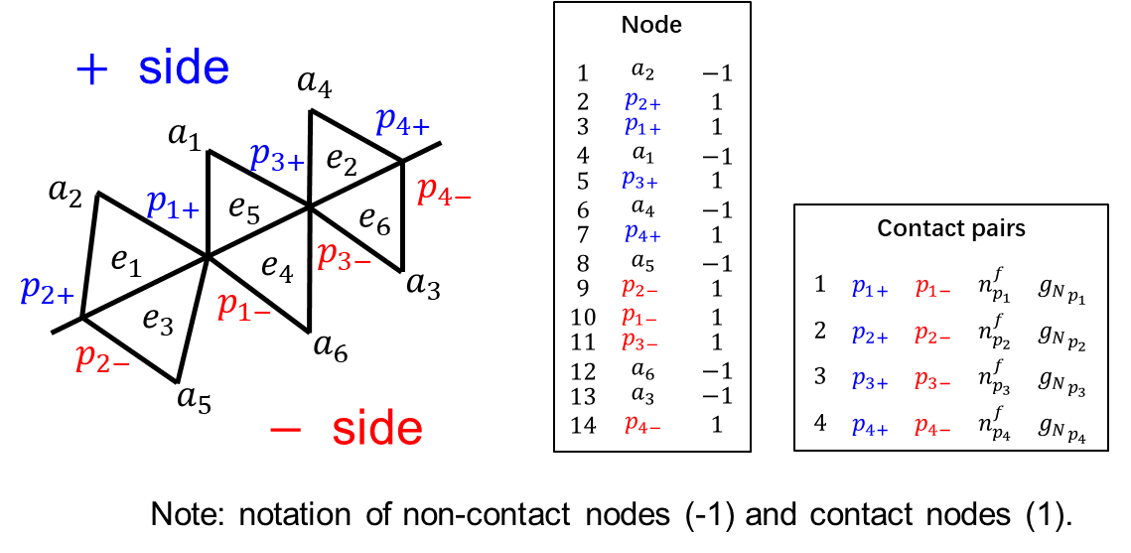}
\caption{ Definition of contact pairs in a local system on the unstructured grids. }
\label{fig:fig_ContactPairsSingleFrac}
\end{figure}

The relative displacement between $+$ side and $-$ side is measured by the jump displacement vector $[\![ {\bm{u}}^f ]\!]$ in local system, as shown in Fig. \ref{fig:fig_ContactPairsSingleFrac}. 
The transformation matrix $\bm{G}$ achieves the expression of relative displacement on contact surface. To convert the quantity from global system to local system, the rotation matrix $\bm{S}$ is defined by unit vectors of the local system. Consequently, the jump displacement reads: 
\begin{equation}
[\![ \hat{\bm{u}}^f ]\!] = \bm{S}^T \bm{G} \bm{U}^f
\label{eq:LocalDisp}
\end{equation}
with the displacement vector $\bm{U}^f$ related to the fracture surface. 
Strain $\bm{\varepsilon}$ and stress $\bm{\sigma}$ are expressed as: 
\begin{equation}
\begin{aligned}
\bm{\varepsilon} = \bm{B} \bm{U} \;\; {\rm and} \;\;
\bm{\sigma} = \bm{D} \bm{B} \bm{U} \\
\end{aligned}
\label{eq:StressStrain}
\end{equation}
with the elastic matrix $\bm{D}$ and the strain operator $\bm{B}$.

It is straightforward that the discretized forms can be obtained once substituting the quantities in Eqs. (\ref{eq:WeakFormState}), (\ref{eq:WeakFormSlip_Disp}) and (\ref{eq:WeakFormSlip_Lambda}) using the approximations Eqs. (\ref{eq:PrimaryApproximations}), (\ref{eq:LocalDisp}) and (\ref{eq:StressStrain}). 
For the stick state: 
\begin{small}
\begin{equation}
\begin{aligned}
\int_{\Omega} \delta \bm{U}^T \bm{B^T} \bm{D} \bm{B} \bm{U} d\Omega 
- \int_{\Omega} \delta \bm{U}^T \left[\bm{N}^u\right]^T \bm{f} d\Omega
- \int_{\Gamma_{ex}^N} \delta \bm{U}^T \left[\bm{N}^u\right]^T \bar{\bm{t}}^{ex} d\Gamma
\\+ 
\sum_{i=1}^{N^f} \int_{\Gamma_{in,i}^{stick}} \delta  \bm{U}^T \bm{G}^T \bm{S} \bm{N}^{\lambda} \bm{\Lambda} d\Gamma
= 0 
\\
\sum_{i=1}^{N^f} \int_{\Gamma_{in,i}^{stick}} \delta  \bm{\Lambda}^T  \left[ \bm{N}^{\lambda} \right]^T \bm{S}^T \bm{G} \bm{U} d\Gamma 
= - \sum_{i=1}^{N^f} \int_{\Gamma_{in,i}^{stick}} \delta  \bm{\Lambda}^T \left[ \bm{N}^{\lambda} \right]^T \bm{g} d\Gamma 
\end{aligned}
\label{eq:FullyDiscrtizedStick}
\end{equation}
\end{small}

For the slip state, the tangential component which related to $\bm{\Lambda}$ in Eq. (\ref{eq:FullyDiscrtizedStick}) need to be replaced by the critical traction $\tau_c$ defined in Eq. (\ref{eq:KKTCondition_tangential}). Meanwhile, it is important that employing Eq. (\ref{eq:IndicatorFunction}) to determine the direction of slip contact.

Eq. (\ref{eq:FullyDiscrtizedStick}) for stick contact as well as its modified version for slip contact are the core for solution strategy of frictional contact and sliding in Section \ref{section:SecSolutionStrategy}. Before that, the treatment of contact and sliding on crossing fractures should be addressed in Section \ref{subsection:MultiCrossingFrac}.


\subsection{A novel treatment of contact on crossing fractures}
\label{subsection:MultiCrossingFrac}

Normally, in classical computational contact mechanics, the concept of contact pair is used to capture the contact or impact between two separate deformable (or rigid) bodies, in which the situation of crossing contact surfaces is unusual, see \cite{Wriggers2006,Zhong1993,Meguid2008}. In contrast to that, it is very common in geoscience application that discrete fractures would cross each other, so that the contact mechanics of crossing fractures should be addressed.

\begin{figure}[H]
\centering
\includegraphics[width=6.5cm]{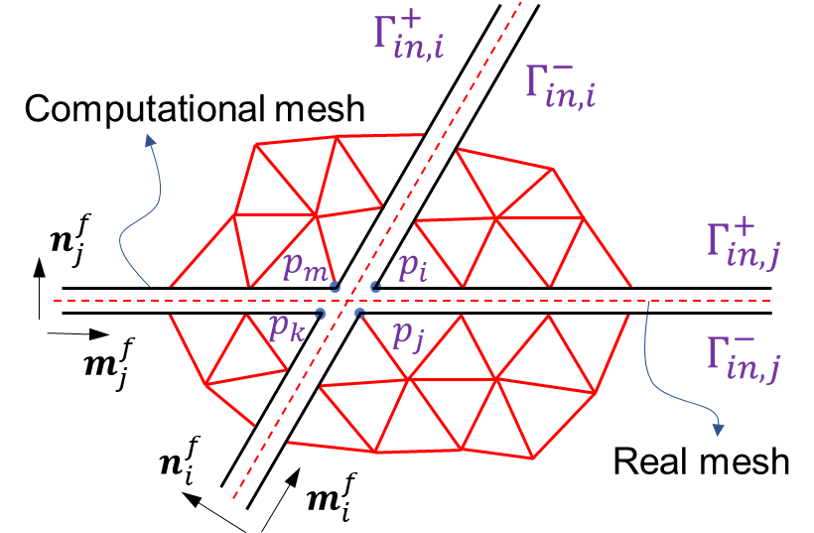}
\caption{ Schematic of contact pairs on crossing fractures. }
\label{fig:fig_CrossingFrac}
\end{figure}

Fig. \ref{fig:fig_CrossingFrac} shows a schematic of crossing fractures and their contact pairs defined in the local systems of fractures $\Gamma_{in,i}$ and $\Gamma_{in,j}$. 
The node-split technique \cite{Aagaard2013} is applied to duplicate the nodes on contact surface, so that the computational mesh allows a relative displacement on the position of crossing fractures. For example, four points $p_i,p_j,p_k,p_m$ are defined on the crossing position. They would be contacted if fracture $\Gamma_{in,i}$ (or fracture $\Gamma_{in,j}$) is imposed by compression. To capture this interaction, the "crossing contact pairs" are devised in the presented numerical scheme. Assuming contact surface of $\Gamma_{in,i}$ would be contacted together, we designed the crossing contact pairs $p_i - p_m $ and $p_k - p_j$. If compression loading is not active, the crossing contact pairs allow to be separated.


\section{Solution strategy of the mixed-FEM scheme}
\label{section:SecSolutionStrategy}

In this section, the algebraic form of the contact system is derived and resulting an unified matrix formulation. 
The system of contact mechanics is solved by iteration method. The preconditioner is introduced to handle the saddle-point algebraic system. It leads to a preconditioned mixed-FE scheme. All algorithms are implemented in our C++ code.


\subsection{The iteration method for contact mechanics}
\label{subsection:SubSecIteration}

The system of contact mechanics is a nonlinear system, where the unknowns $\bm{U}$ and $\bm{\Lambda}$ are coupled together. 
As shown in system Eq. (\ref{eq:FullyDiscrtizedStick}), the two equations depend on each other. Consequently, the mixed-finite element formulation \cite{Cescotto1993,Auricchio2017,Franceschini2019} is formulated. In this way, the two primary unknowns can be calculated by one non-linear system in a monolithic strategy. 
Residual vector $\bm{R} = [\bm{R}^u \;\; \bm{R}^{\lambda}]^T$ is constructed for the iteration. The superscripts $u$ and $\lambda$ mean the displacement and Lagrange multiplier. In the iteration, we wish the $L^2$-norm of residual vector converges to zero: 
\begin{equation}
\lim_{n \rightarrow \infty} \Vert \bm{R}^{n} \Vert_2 = \bm{0}
\end{equation}
where the components of residual vector $\bm{R}$ are defined based on Eq. (\ref{eq:FullyDiscrtizedStick}): 
\begin{small}
\begin{equation}
\begin{aligned}
\bm{R}^u = 
\int_{\Omega} \delta \bm{U}^T \bm{B^T} \bm{D} \bm{B} \bm{U} d\Omega 
- \int_{\Omega} \delta \bm{U}^T \left[\bm{N}^u\right]^T \bm{f} d\Omega
- \int_{\Gamma_{ex}^N} \delta \bm{U}^T \left[\bm{N}^u\right]^T \bar{\bm{t}}^{ex} d\Gamma
\\- 
\sum_{i=1}^{N^f} \int_{\Gamma_{in,i}^{stick}} \delta  \bm{U}^T \bm{G}^T \bm{S} \bm{N}^{\lambda} \bm{\Lambda} d\Gamma
\\
\bm{R}^{\lambda}  =
\sum_{i=1}^{N^f} \int_{\Gamma_{in,i}^{stick}} \delta  \bm{\Lambda}^T  \left[ \bm{N}^{\lambda} \right]^T \bm{S}^T \bm{G} \bm{U} d\Gamma 
+ \sum_{i=1}^{N^f} \int_{\Gamma_{in,i}^{stick}} \delta  \bm{\Lambda}^T \left[ \bm{N}^{\lambda} \right]^T \bm{g} d\Gamma 
\end{aligned}
\label{eq:ResidualStick}
\end{equation}
\end{small}

Then, the Newton-Raphson iteration method is introduced to linearize the nonlinear system. Jacobian $\bm{J}^{\nu}$ and residual vector $\bm{R}^{\nu}$ should be calculated at each iteration step $\nu$. The unknown vector $\bm{x} = [\bm{U} \;\; \bm{\Lambda} ]^T$ can be solved at step $\nu+1$ in incremental form $\delta \bm{x}^{\nu+1}$: 
\begin{equation}
\begin{bmatrix}
\bm{J}_{11} & \bm{J}_{12}  \\
\bm{J}_{21} & \bm{J}_{22} 
\end{bmatrix}^{\nu}
\begin{bmatrix}
\delta \bm{U} \\
\delta \bm{\Lambda} \\
\end{bmatrix}^{\nu+1}
= - 
\begin{bmatrix}
\bm{R}^u \\
\bm{R}^{\lambda} \\
\end{bmatrix}^{\nu}
\label{eq:SystemMatrixForm}
\end{equation}
with the components of Jacobian $\bm{J}_{ij}$ $\left( i,j=1,2 \right)$, which are determined by the derivative of residual vector with respect to unknown vector: 
\begin{equation}
\bm{J}^{\nu} = \left. \frac{\partial \bm{R}}{\partial \bm{x}} \right \vert ^{\nu}
\label{eq:Jacobian}
\end{equation}

The contact system can be then written as the algebraic form based on Eqs. (\ref{eq:SystemMatrixForm}) and (\ref{eq:Jacobian}): 
\begin{equation}
\begin{bmatrix}
\bm{K}^{uu} & \bm{C}^{\lambda u}  \\
\left[\bm{C}^{\lambda u}\right]^T & \bm{0} 
\end{bmatrix}
\begin{bmatrix}
\delta \bm{U} \\
\delta \bm{\Lambda} \\
\end{bmatrix}
= - 
\begin{bmatrix}
\bm{R}^u \\
\bm{R}^{\lambda} \\
\end{bmatrix}
\label{eq:SaddlePointSystem}
\end{equation}
where the effects induced by slip and stick could be reflected through block $\bm{C}^{\lambda u}$ and its transpose $\left[\bm{C}^{\lambda u}\right]^T$. $\bm{K}^{uu}$ is the block related to displacement. 
The components in Eq. (\ref{eq:SaddlePointSystem}) are provided in \ref{appendix_discretizedForms}.

\subsection{Preconditioning of the saddle-point algebraic system}
\label{subsection:SubsecPreconditioning}

\begin{figure}[H]
\centering
\includegraphics[width=7cm]{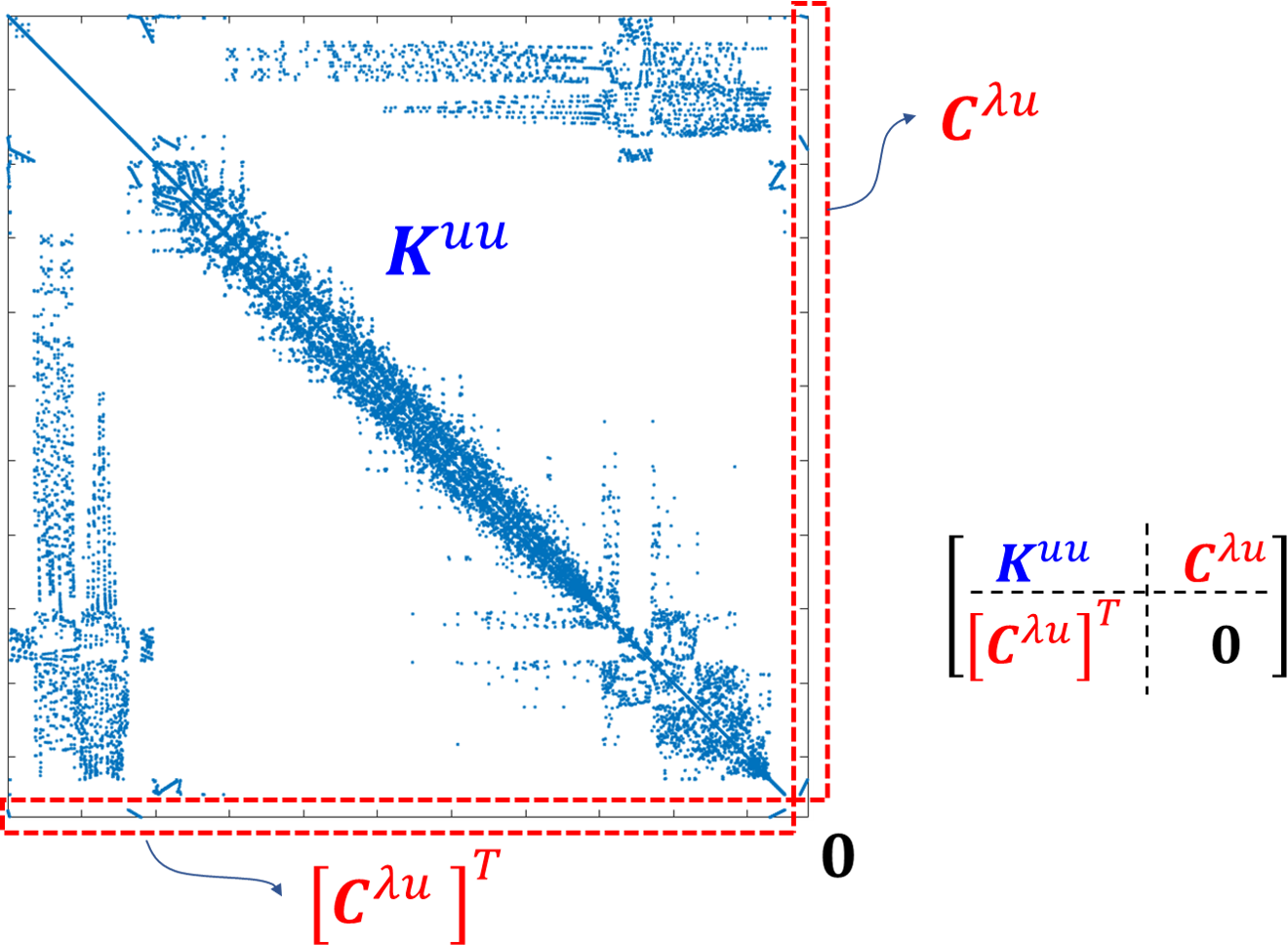}
\caption{ The pattern of non-zero entries of saddle-point algebraic system in the presented numerical scheme. }
\label{fig:fig_SaddlePointStructure}
\end{figure}

It should be noted that the lower block diagonal entry of the Jacobian in Eq. (\ref{eq:SaddlePointSystem}) is a zero matrix.
The shape of $\bm{C}^{\lambda u}$ is rectangular, while $\bm{K}^{uu}$ is a square matrix, as shown in Fig. \ref{fig:fig_SaddlePointStructure}. 
These features lead to a special algebraic structure, which is the so-called saddle-point system \cite{Franceschini2019,Auricchio2017,Pestana2015}. 
Furthermore, the ill-condition Jacobian $\bm{J}$, which has a high condition number, would lead to a numerical instability when solving the saddle-point system Eq. (\ref{eq:SaddlePointSystem}). 
To this end, the precondition technique is presented to improve the numerical robustness of the system.

The preconditioned Jacobian $\bar{\bm{J}}$ with a low condition number is derived through the preconditioned operation: 
\begin{equation}
\bar{\bm{J}} = \bm{P} \bm{J} = 
\begin{bmatrix}
\bm{A}^{-1} & \bm{0} \\
\bm{0} & \bm{B}^{-1}
\end{bmatrix}
\begin{bmatrix}
\bm{K}^{uu} & \bm{C}^{\lambda u}  \\
\left[\bm{C}^{\lambda u}\right]^T & \bm{0} 
\end{bmatrix}
\label{eq:PreconditionedJacobian}
\end{equation}
with the components $\bm{A}$ and $\bm{B}$ of preconditioner $\bm{P}$. They are defined by the norms of $\bm{K}^{uu}$ and $\bm{C}^{\lambda u}$. For convenience, we denote $\left[\bm{C}^{\lambda u}\right]^T = \bm{D}$ and $\left[ \bm{K}^{uu} \;\; \bm{C}^{\lambda u} \right] = \bm{E}$. We refer to \ref{appendix_preconditioner} for the expanded forms. 

Combining the system of Eqs. (\ref{eq:SaddlePointSystem}), (\ref{eq:Blocks1}), (\ref{eq:Blocks2}) with the operation Eq. (\ref{eq:PreconditionedJacobian}), a preconditioned mixed-FE scheme is presented, which is one of the core innovations of the presented work. The attractive feature is ease of implementation and it can be directly integrated into an existing solver.

\subsection{Global strategy: the monolithic-updated contact algorithm}
\label{subsection:GlobalAlgorithm}

The contact system should be resolved using an iteration strategy, since the governing equations of contact mechanics are a set of nonlinear partial differential equations. Moreover, the contact state of each contact pair should be updated dynamically at each time step and depend on the loading condition. It is essential that to check the contact state based on the contact constrained conditions and current stress/displacement condition. Meanwhile, each of the contact pairs should be evaluated at each iteration. To this end, the monolithic-updated contact algorithm is designed to update the unknowns (displacement and Lagrange multiplier) in one algebraic system, as shown in \ref{appendix_algorithm}.


\section{Numerical results and discussion}
\label{section:NumericalTests}

In this section, a series of numerical tests is conducted to study the frictional contact and sliding of single- and multi-crossing fractures based on the proposed mixed-FE scheme. First, a benchmark study is presented to verify the validation. Later, the contact behavior on single fracture is studied under mixed mode loading. Finally, tests with complex geometry are studied. The slippage and opening are analyzed with different conditions.


\subsection{Fracture slip controlled by frictional law}
\label{subsection:Benchmark1}

\begin{figure}[H]
\centering
\includegraphics[width=11cm]{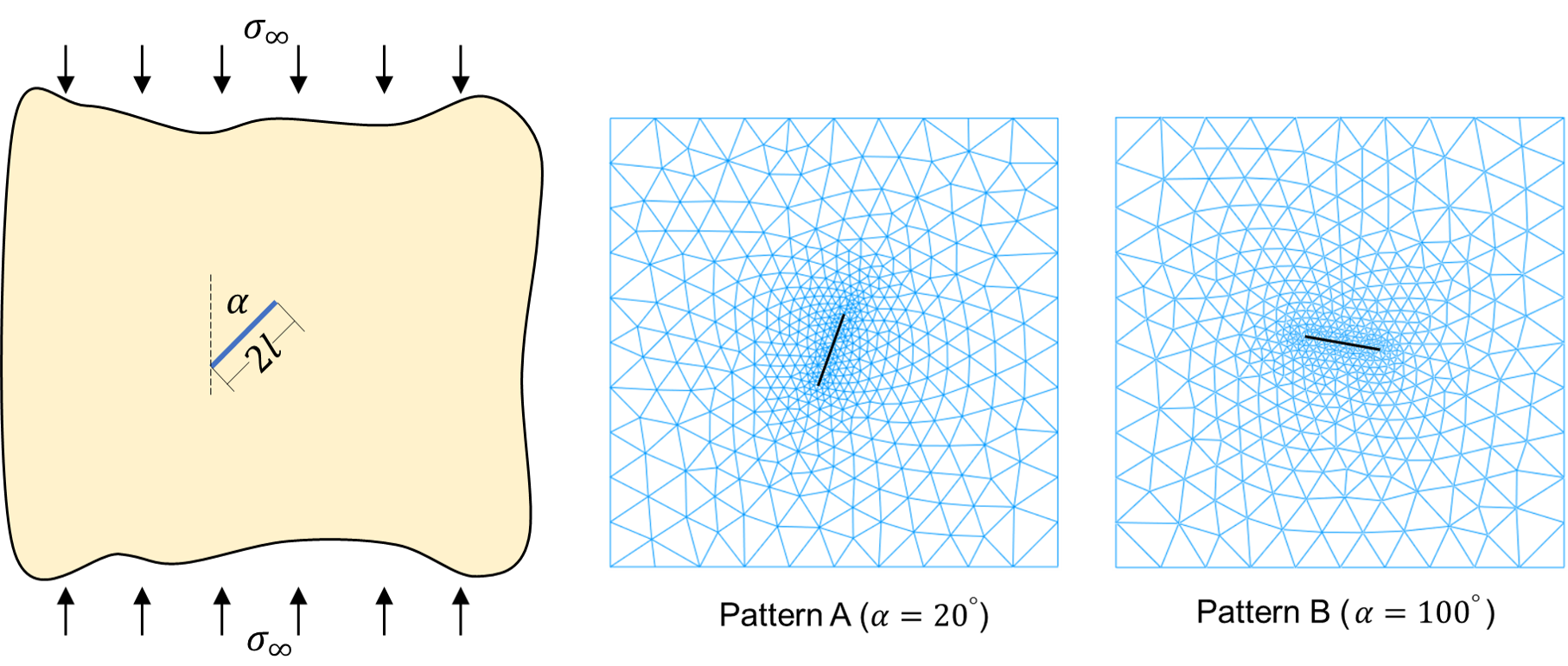}
\caption{ Schematic of the model with different angles. The contact behavior is controlled by frictional law. }
\label{fig:fig_Benchmark1_ModelMesh}
\end{figure}

A fractured medium, which contains a single inclined fracture with angle $\alpha$ and length $2l$, is modeled when the uniaxial loading $\sigma_{\infty}$ is imposed. As shown in Fig. \ref{fig:fig_Benchmark1_ModelMesh}, the numerical model is constructed and two patterns with different inclined angles are considered. We adopt the assumptions of homogeneous and linear elastic material with Young's modulus $E = 25 {\rm GPa}$ and Poisson ratio $\nu = 0.25$. The parameters of frictional law are frictional angle $\varphi = 30^{\circ}$ and cohesion $c=0$. Note that the compression $\sigma_{\infty} = 10 {\rm MPa}$ and crack length $2l = 2 {\rm m}$.


This model is taken as a benchmark test to verify the proposed numerical method. To this end, the slippage (i.e. the relative tangential displacement) on contact surface can be calculated by the analytical solution from literature \cite{Phan2003}: 
\begin{equation}
[\![ {u}^f_T ]\!] = \frac{4 t_T^f \left( 1-\nu^2 \right)}{E} \sqrt{l^2 - \left( \eta-l \right)^2}
\label{eq:Benchmark1_Disp}
\end{equation}
with the tangential and normal components of contact traction: 
\begin{equation}
\begin{aligned}
&t_T^f = \sigma_{\infty} \sin\alpha \cos\alpha - \sigma_{\infty} \sin^2 \alpha \tan\varphi \\
&t_N^f = - \sigma_{\infty} \sin^2{\alpha} 
\end{aligned}
\label{eq:Benchmark1_Traction}
\end{equation}
where $ 0 \leqslant \eta \leqslant 2l$ is the coordinate along the length of contact surface.

\begin{figure}[H]
\centering
\includegraphics[width=12cm]{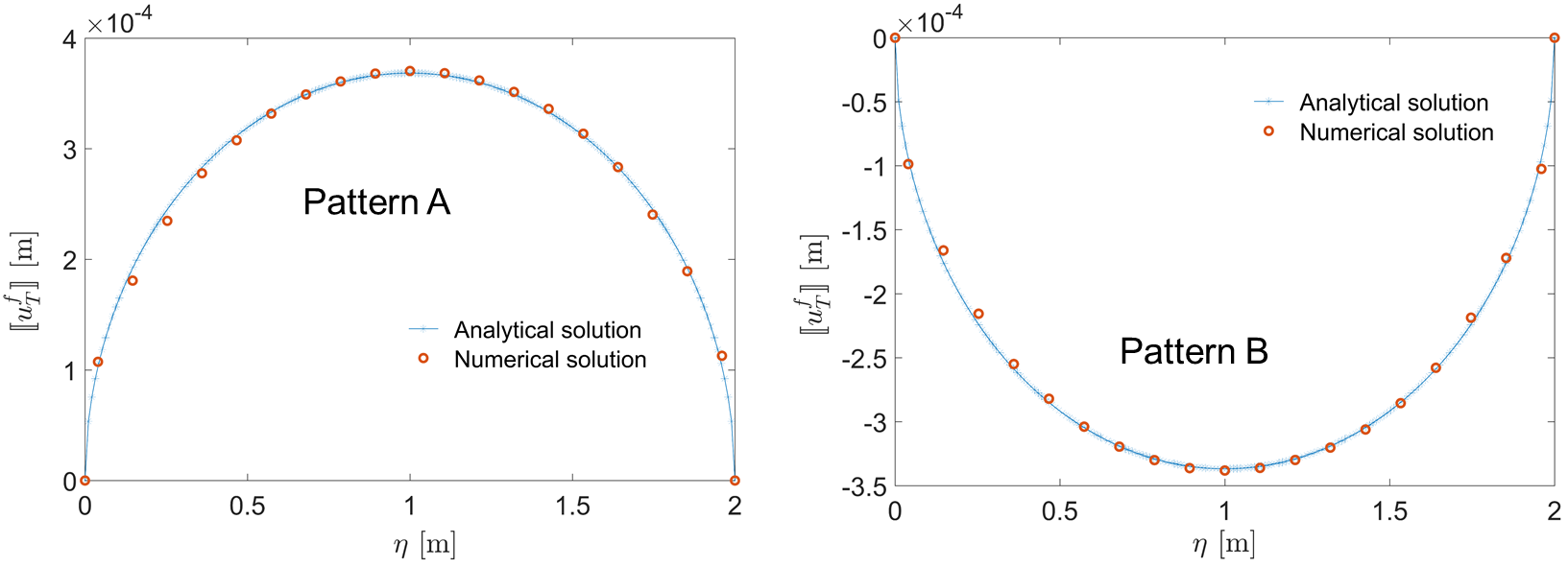}
\caption{ The comparison between analytical and numerical solutions of single inclined fracture. Slip on fracture for Pattern A and Pattern B. }
\label{fig:fig_benchmark1_slip}
\end{figure}

\begin{figure}[H]
\centering
\includegraphics[width=12cm]{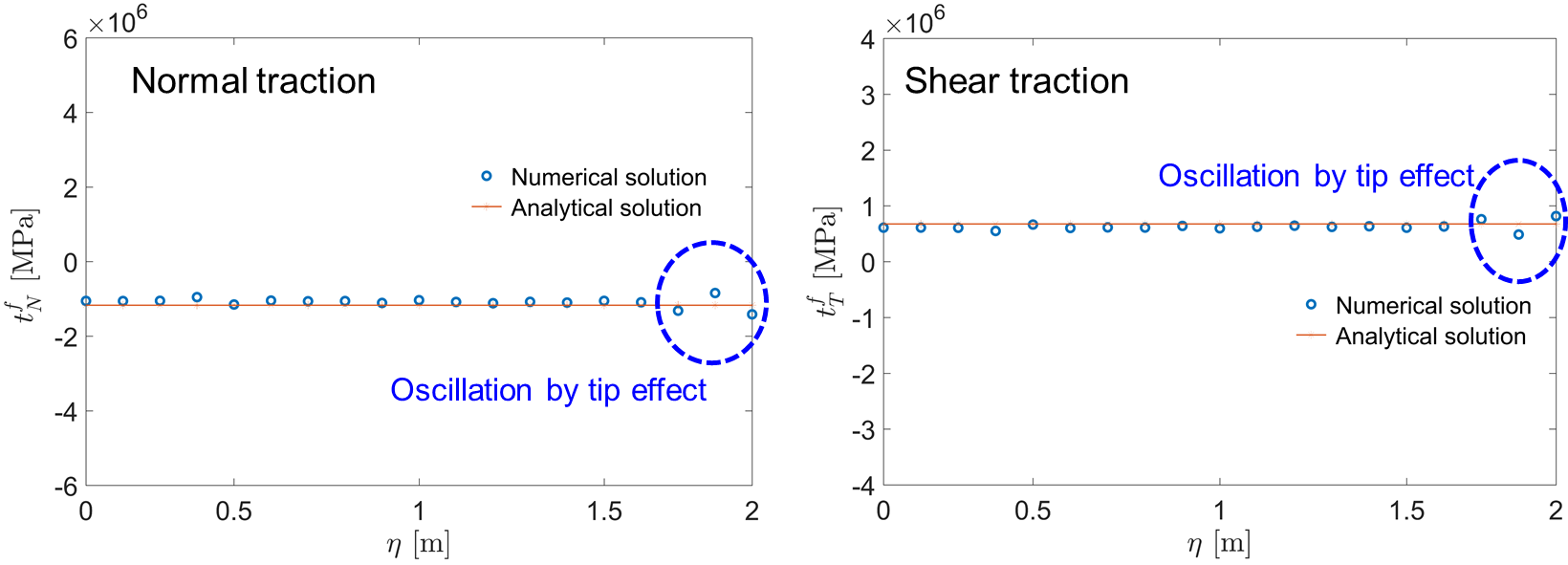}
\caption{ The comparison between analytical and numerical solutions of the single inclined fracture. Contact traction on fracture for Pattern A. }
\label{fig:fig_benchmark1_traction}
\end{figure}

Fig. \ref{fig:fig_benchmark1_slip} (left) illustrates the comparison between the analytical and numerical solutions. As shown in this figure, the relative displacement (slippage) on contact surface shows a "parabolic" shape. Fig. \ref{fig:fig_benchmark1_traction} indicates the components of traction maintain constant. The analytical and numerical results agree well with each other. Moreover, the oscillation on the crack tip is observed in the endpoint of x-coordinate, which also has been reported in literature \cite{Franceschini2016}. 
Furthermore, the tangential component of relative displacement on contact surface is calculated in case of Pattern B, as shown in Fig. \ref{fig:fig_benchmark1_slip} (right). The results imply that the numerical solution is consistent with the analytical results. 

\begin{figure}[H]
\centering
\includegraphics[width=8cm]{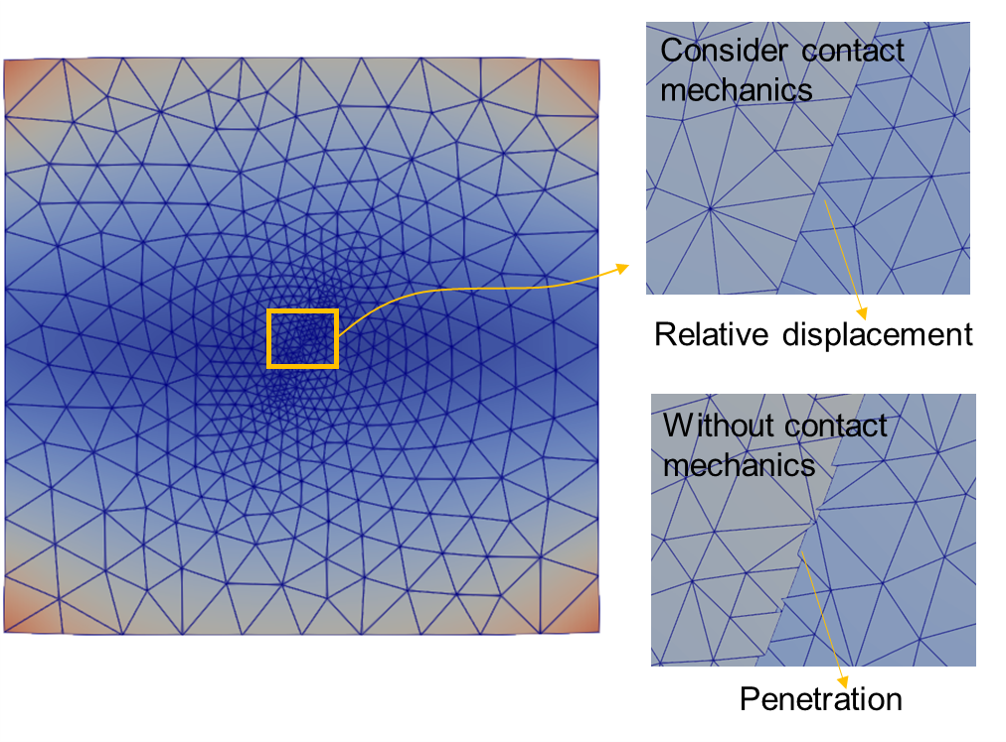}
\caption{ Deformation of finite element grids in different situations. The numerical treatment of frictional contact and sliding prevents the penetration of grid on two sides of the contact surface. }
\label{fig:Enlargedfig_Benchmark1_PatternA_ZOOMView}
\end{figure}

The importance of considering contact behavior is depicted in Fig. \ref{fig:Enlargedfig_Benchmark1_PatternA_ZOOMView}. As illustrated in the enlarged view, the grids on two sides of contact surface would be penetrated to each other if contact constraints are neglected. In contrast, a pair of contact traction corresponding to the two sides would be naturally created in the framework of Lagrange multiplier method to prevent the penetration when considering contact effect on fracture.

\subsection{Shear failure on contact surface}

A fractured medium, which is intersected by a single fracture crossing the entire domain, is used to analyze the pure slippage along the contact surface. We follow the same parameters as the benchmark in literature \cite{Borja2008}. The parameters of elastic material are $E = 5 $ GPa and $\nu = 0.3$. The frictional angle $\varphi = 5.71^{\circ}$ (means the frictional coefficient $\tan\varphi \approx 0.1$) and cohesion $c=0$.

\begin{figure}[H]
\centering
\includegraphics[width=10.5cm]{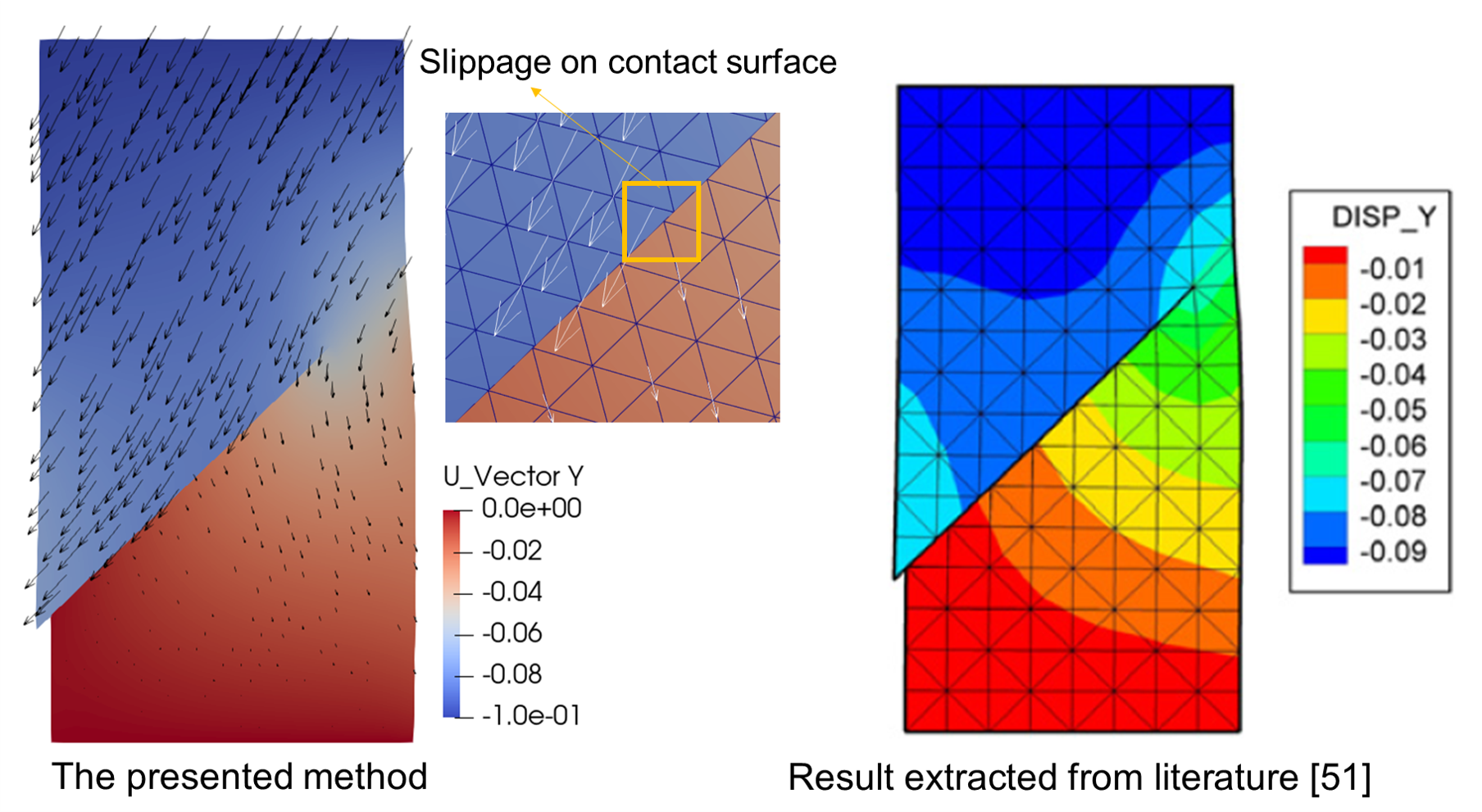}
\caption{ Displacement pattern of shear failure model (Pattern A). The slippage of contact surface which cuts the entire domain is measured by the tangential displacement.   }
\label{fig:fig_Benchmark3_Umagnitude}
\end{figure}

\begin{figure}[H]
\centering
\includegraphics[width=11.5cm]{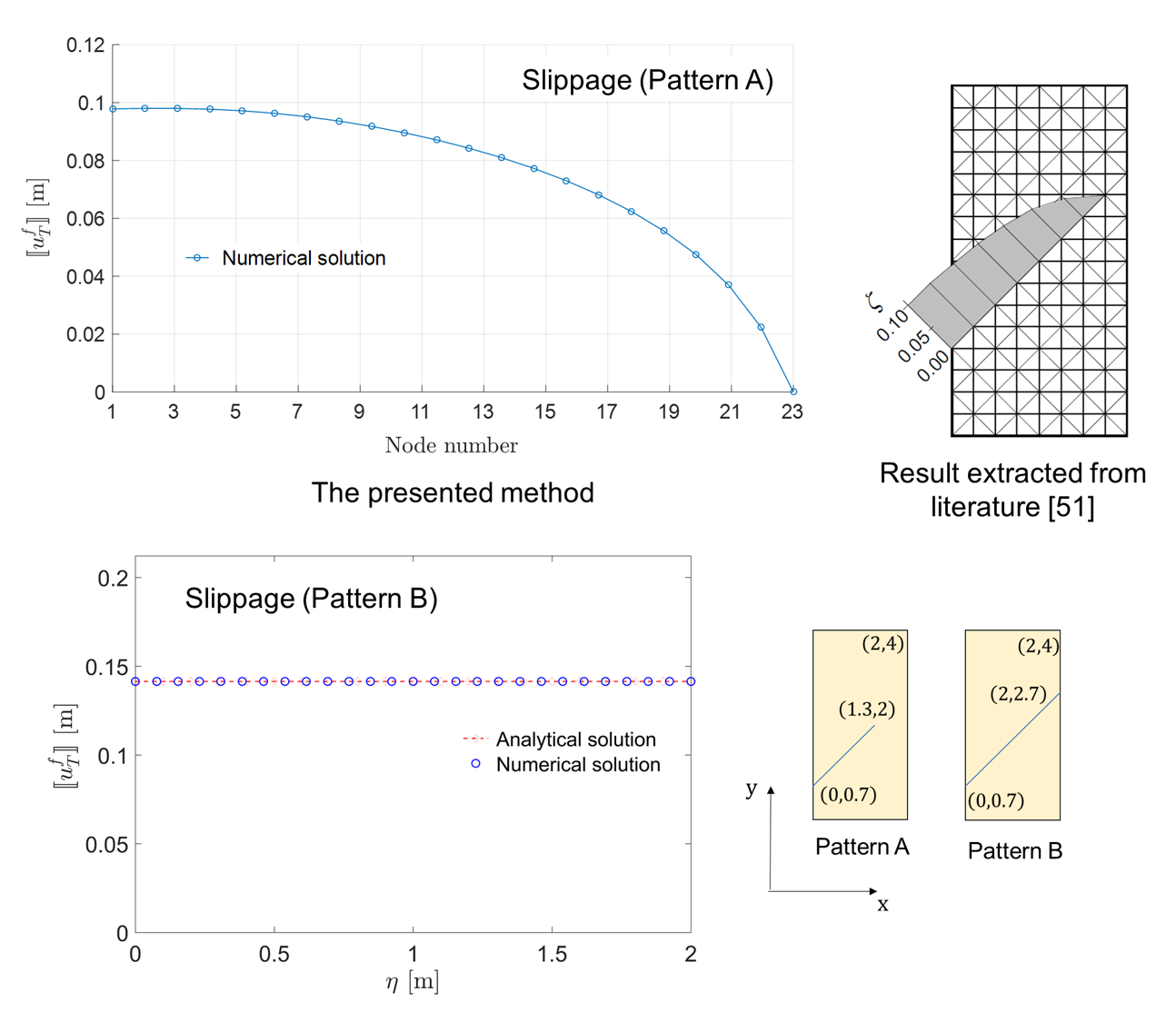}
\caption{ The comparison between numerical solution and reference solution of Pattern A (top) and Pattern B (bottom). }
\label{fig:fig_Benchmark3_curves_uT}
\end{figure}

In this situation (Pattern A), the crack tips are no longer surrounded by the host elastic matrix, thus deformation contour shows a different pattern compared with the case of embedded-fracture model (for instance the case in Section \ref{subsection:Benchmark1}), as displayed in Fig. \ref{fig:fig_Benchmark3_Umagnitude}. The upper part of the geometry has a significant movement, such that the entire body slips along the fracture surface. Fig. \ref{fig:fig_Benchmark3_curves_uT} shows a comparison of results calculated by the presented method and the reference solution.

Another setting of this model is presented in \cite{Franceschini2020}, denoted by Pattern B in this section. The analytical solution of slippage is a constant $\Delta u \approx 0.1414 {\rm m}$ according to the literature \cite{Borja2008,Franceschini2020}. Fig. \ref{fig:fig_Benchmark3_curves_uT} illustrates that the analytical and numerical solutions agree with well each other.

\subsection{The effect of fluid pressure and mixed mode loadings}

Different types of loading mode would lead to various responses of the fractured medium. The preceding case studies focus on pure uniaxial compression or single type of loading instead of complex loading condition. In contrast to that, the compressive+shear mode (Pattern A) and the tensile+shear mode (Pattern B) are simulated. $E = 50$ GPa and $\varphi = 30^{\circ}$. The normal stress equals $10$ MPa on the top surface of the model, while the shear stress allows varied $5 \sim 8$ MPa. It is obvious that the contact constraints would be active if Pattern A is applied. While the situation shows a distinct in Pattern B, in which the two sides of fracture surface would be separate thus the contact constraints are no longer active.

\begin{figure}[H]
\centering
\includegraphics[width=13.5cm]{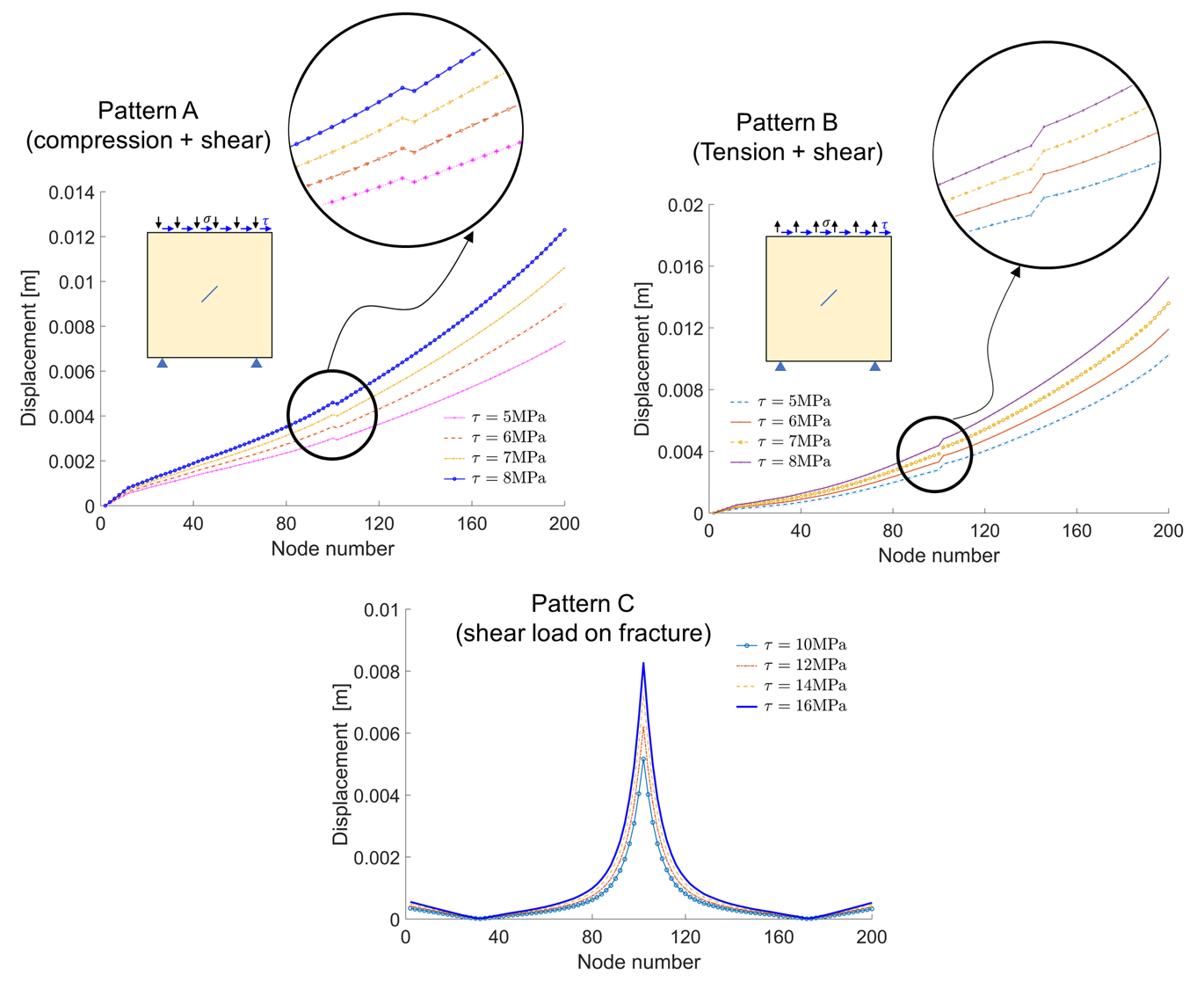}
\caption{ The variation of displacement across the contact surface in Pattern A (left), Pattern B (right) and Pattern C (bottom). }
\label{fig:fig_Test4_PatternABC_curve}
\end{figure}

Fig. \ref{fig:fig_Test4_PatternABC_curve} (top) provides several key features and gives a comparison between Pattern A and Pattern B. Note that the jump displacement at the position which is intersected by the fracture exhibits different behaviors as shown in the enlarged views. Furthermore, the increased value of shear loading leads to the corresponding increased displacement magnitude on fracture.

\begin{figure}[H]
\centering
\includegraphics[width=12.5cm]{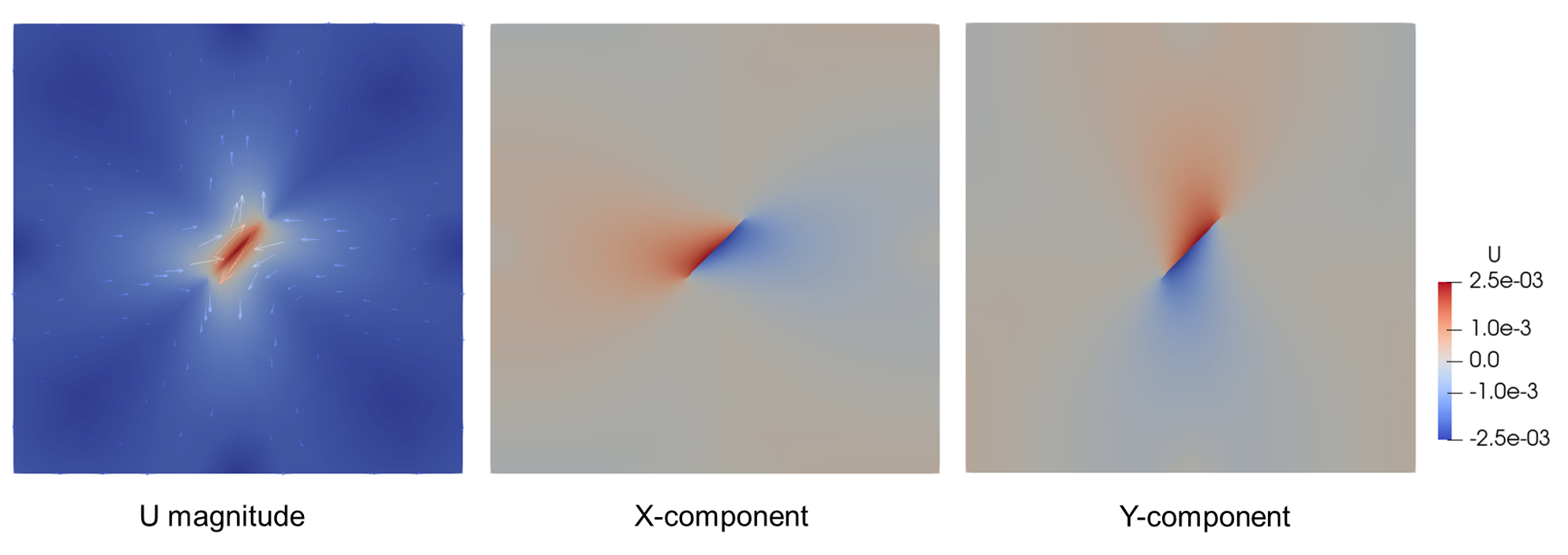}
\caption{ Profiles of Pattern C. The traction on contact surface leads to the symmetric mode of displacement. }
\label{fig:fig_Test4_PatternC_countor}
\end{figure}

Pattern C is designed to analyze the shear behavior on contact surface. Fig. \ref{fig:fig_Test4_PatternC_countor} shows the deformation profile when fracture surface is imposed by an inverse pair of traction on the corresponding two sides of fracture surface. The range of shear stress is $10 \sim 16$ MPa. It can be seen from this figure that the shear loading produces a symmetric distribution of displacement field along the diagonal line in the domain. A monitoring line is placed along the off-diagonal line in the domain to measure the variation of displacement. Fig. \ref{fig:fig_Test4_PatternABC_curve} (bottom) captures a sharp shift of displacement magnitude on the fracture.

\begin{figure}[H]
\centering
\includegraphics[width=11cm]{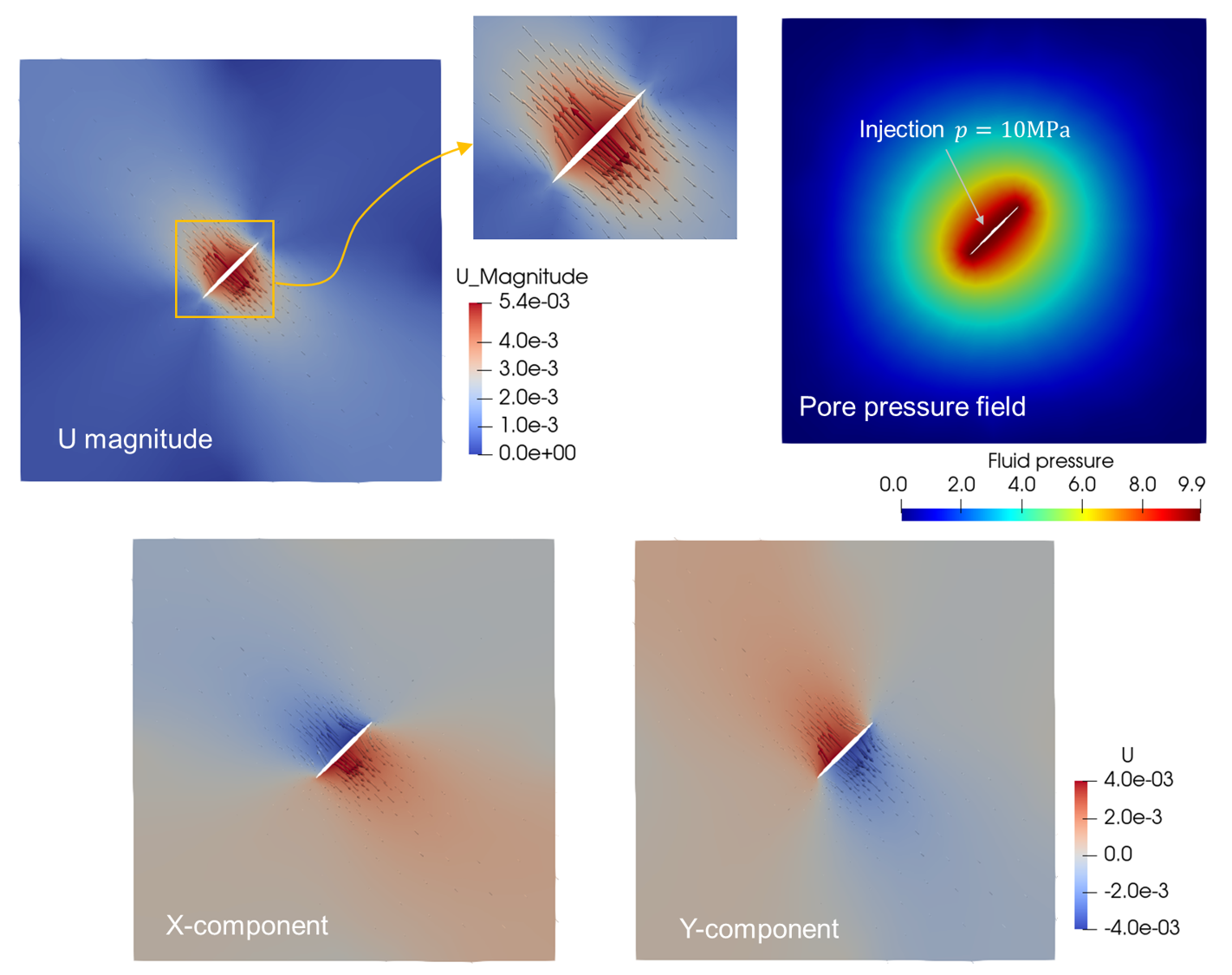}
\caption{ Deformation profiles when fluid pressure applied on fracture plane. Note that the values are amplified by 10 times to show the deformable contours. }
\label{fig:fig_benchmark2_countor}
\end{figure}

The scenario is different when the fluid pressure is considered. As displayed in Fig. \ref{fig:fig_benchmark2_countor}, the fractured surface is imposed by internal pressure which is normal to the two side of fracture, and their directions are opposite to each other. In this way, one can model the scenario of fluid pressure in a hydro-mechanical process. The displacement vector field in Fig. \ref{fig:fig_benchmark2_countor} illustrates the opening of fracture surface. Especially, the discontinuous phenomenon generated by the opened surface can be calculated by an analytical solution from literature \cite{Ucar2018,Sneddon1995}: 
\begin{equation}
\Delta u = \frac{2 l p \left( 1-\nu \right) }{G} \sqrt{ 1- \left( \frac{\eta}{l} \right)^2 }
\label{eq:Benchmark2_Opening}
\end{equation}
with the applied fluid pressure $p = 10$ MPa and coordinate $\eta$ on contact surface. $l = 1$ m is the half length of the fracture, as shown in Fig. \ref{fig:fig_Benchmark1_ModelMesh}. $G$ is the shear modulus of the elastic host matrix. 

\begin{figure}[H]
\centering
\includegraphics[width=7cm]{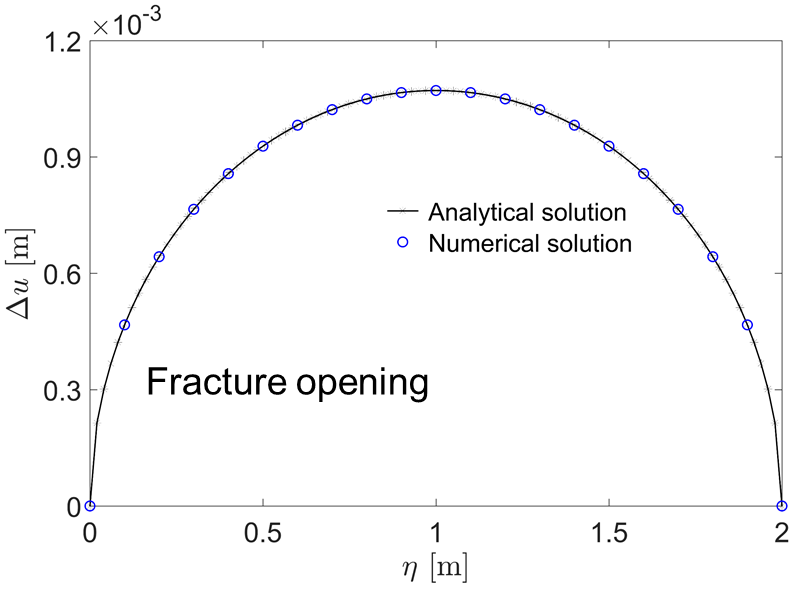}
\caption{ The solution of fracture opening of analytical and numerical solutions when fluid pressure imposed on the fracture. }
\label{fig:fig_Benchmark2_uN_curve}
\end{figure}

Fig. \ref{fig:fig_benchmark2_countor} shows the deformation pattern is symmetric along the off-diagonal line across the entire domain. The components of the up and down sides of the fracture surface illustrate the deformation distribution jumps by the sharp region that intersected by fracture. The enlarged view of the fracture provides an enhanced observation to show the opening and displacement vector field. Fig. \ref{fig:fig_Benchmark2_uN_curve} shows a comparison between analytical and numerical solutions. The shape of opening curve is "parabolic" which is same with the curve of slippage in compressive loading.

\subsection{Crossing fractures under mixed mode loadings}

\begin{figure}[H]
\centering
\subfigure[ Single crossing fractures ]{
\includegraphics[width=3.5cm]{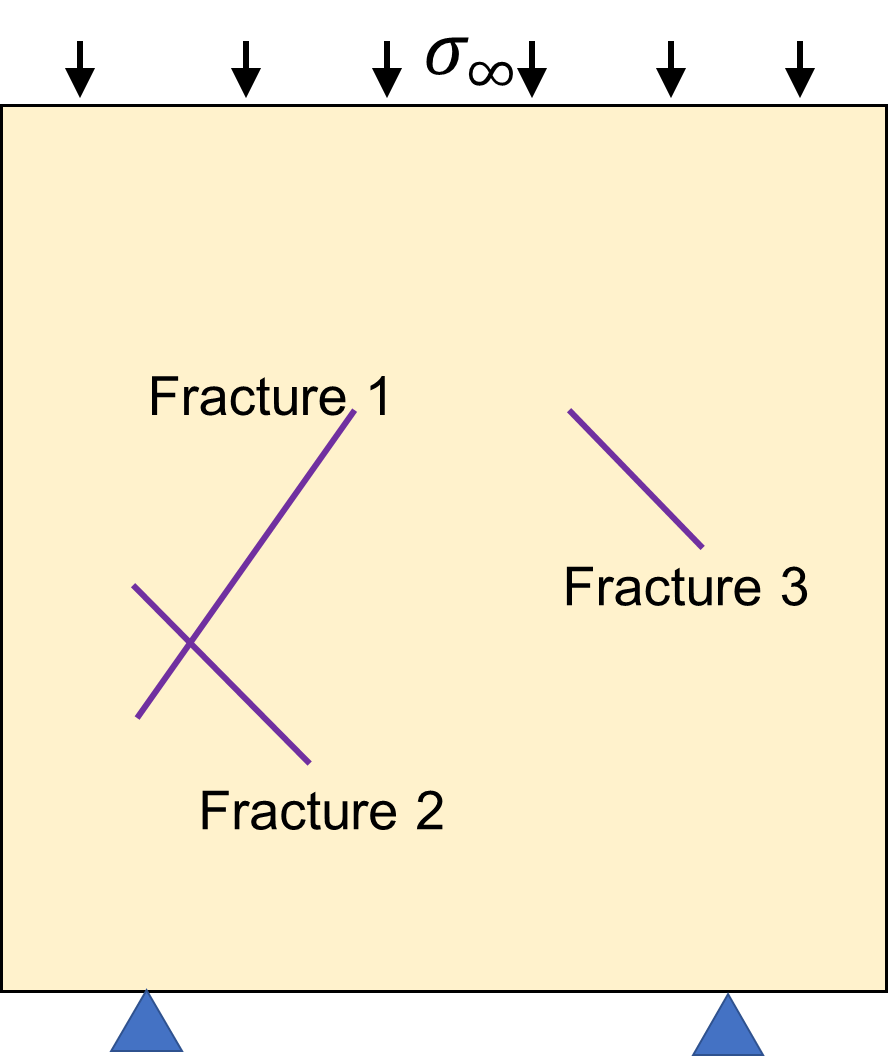}
}
\quad
\subfigure[ Multiple crossing fractures ]{
\includegraphics[width=3.6cm]{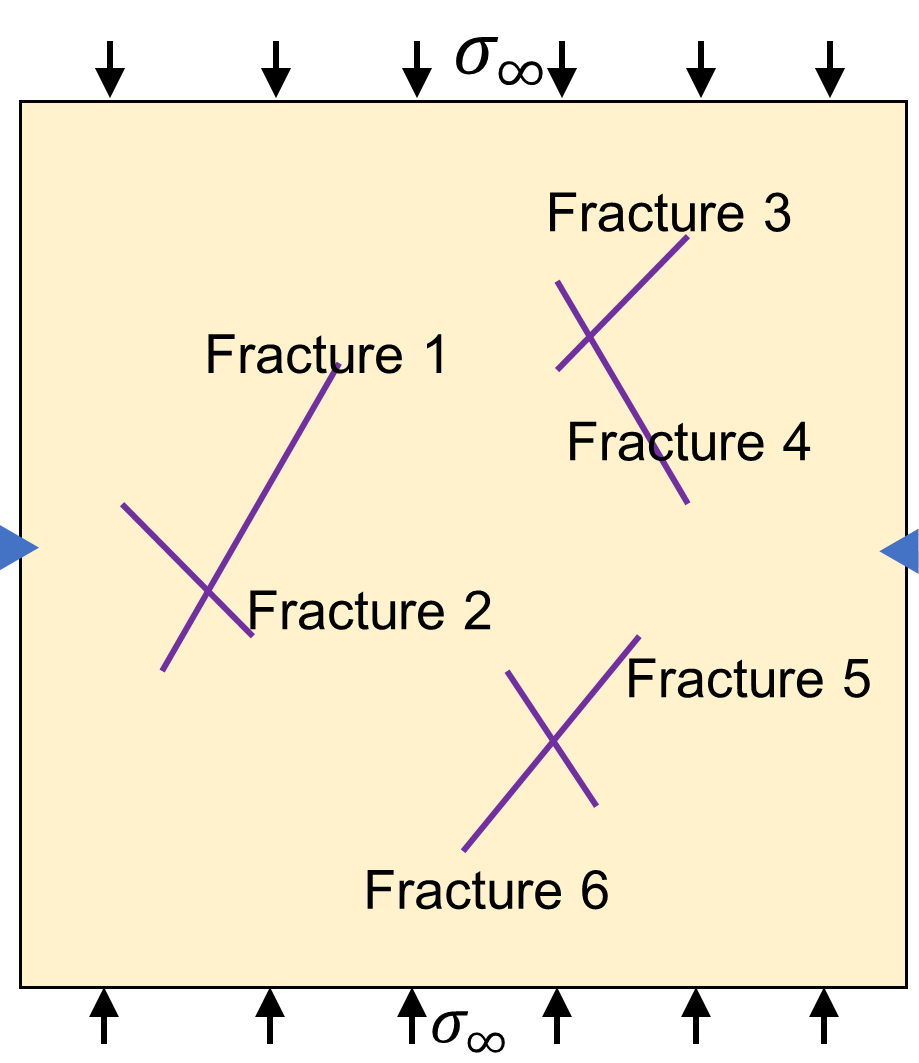}
}
\caption{ Schematics of numerical tests with different fracture patterns. }
\label{fig:fig_Test5and6_GeometricalModel}
\end{figure}

The simulation of frictional contact and shear failure on crossing fractures renders several challenges in computational mechanics \cite{Annavarapu2013,Stefansson2020}. The capability of modeling frictional contact and sliding on multiple crossing fractures is one of the innovations of the proposed work, as presented in Sections \ref{section:SecMathFormulation} and \ref{section:SecDiscretization}. To this end, a complex fractured medium is simulated to show the contact behavior of fractured media containing crossing fractures, as shown in Fig. \ref{fig:fig_Test5and6_GeometricalModel}. The parameters are same as the above case.

\begin{figure}[H]
\centering
\includegraphics[width=15cm]{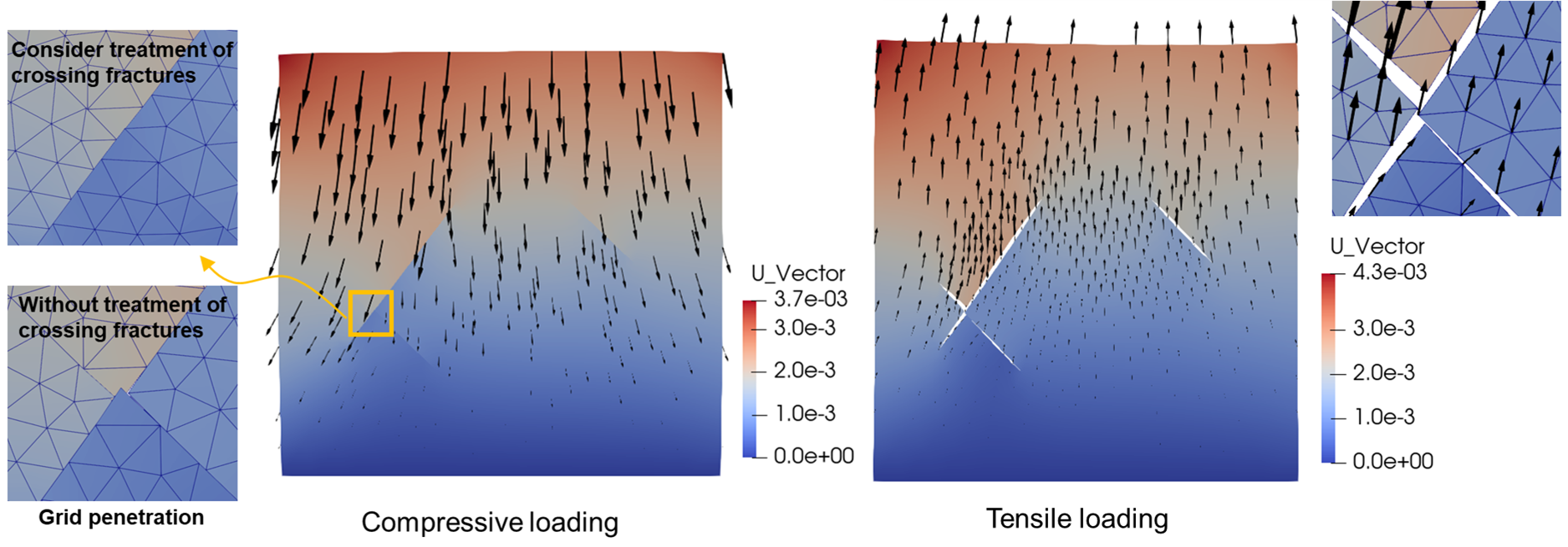}
\caption{ Displacement vector field of the test with single crossing fractures. A comparison of the grid patterns in the left insets shows the effect of numerical treatment on crossing fractures. Note that deformation is amplified by 50 times for visualization. }
\label{fig:fig_crossingFratures_counter}
\end{figure}

The deformation profile, when the compressive loading is applied, is shown in Fig. \ref{fig:fig_crossingFratures_counter} (left). The interaction of multi-fractures at the crossing position can be observed in the enlarged views. A comparison in Fig. \ref{fig:fig_crossingFratures_counter} (the enlarged inset in left) is provided to demonstrate the effect of treatment of crossing fractures on contact mechanics. The displacement vector illustrates the influence of the existing contact surfaces, where a direction shift of the displacement vector is observed around the fractures. Moreover, the treatment of crossing fractures is important to the numerical simulation. Otherwise, a mesh penetration would be occurred at the position of crossing fractures.

\begin{figure}[H]
\centering
\includegraphics[width=15cm]{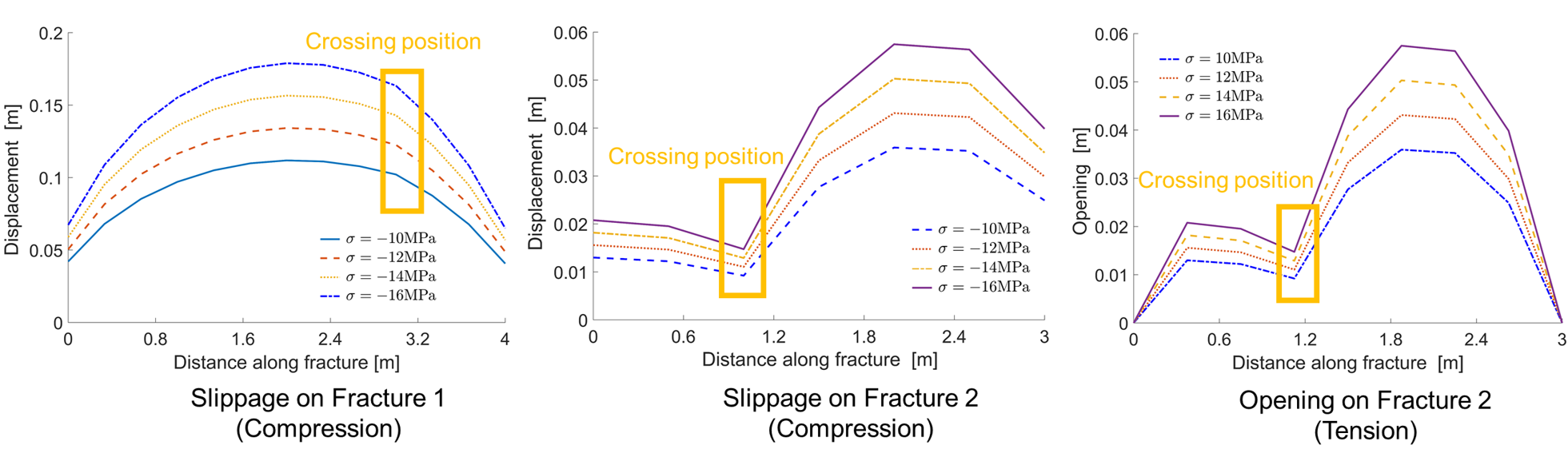}
\caption{  The slippage and opening of fracture surfaces in the single crossing fractures model. }
\label{fig:fig_singleCrossingFrac_disp}
\end{figure}

Fig. \ref{fig:fig_singleCrossingFrac_disp} depicts the variations of slippage on Fractures 1 and 2, which are labeled in Fig. \ref{fig:fig_Test5and6_GeometricalModel}a. When compression, as shown in Fig. \ref{fig:fig_crossingFratures_counter} (left), with the increased stress, the values of slippage would be increased as well. However, the curves of Fractures 1 and 2 show different features, the crossing positions is marked by a box in this figure. It can be observed that the slippage on Fracture 1 is greater than that of Fracture 2. The shape of slippage curve shown in Fig. \ref{fig:fig_singleCrossingFrac_disp} (left) is parabolic and is comparable with the shapes discussed in Section \ref{subsection:Benchmark1}.

As a contrast comparison study, another loading mode is applied to analyze the opening on fracture surface. This distinct can be classified by the vector arrows in the two pictures in Figs. \ref{fig:fig_crossingFratures_counter}. In addition, the opening at the crossing position is observed and shown by the amplified view, which demonstrates a desirable capability of the proposed numerical scheme to handle the opened fracture surfaces.

\begin{figure}[H]
\centering
\includegraphics[width=14.5cm]{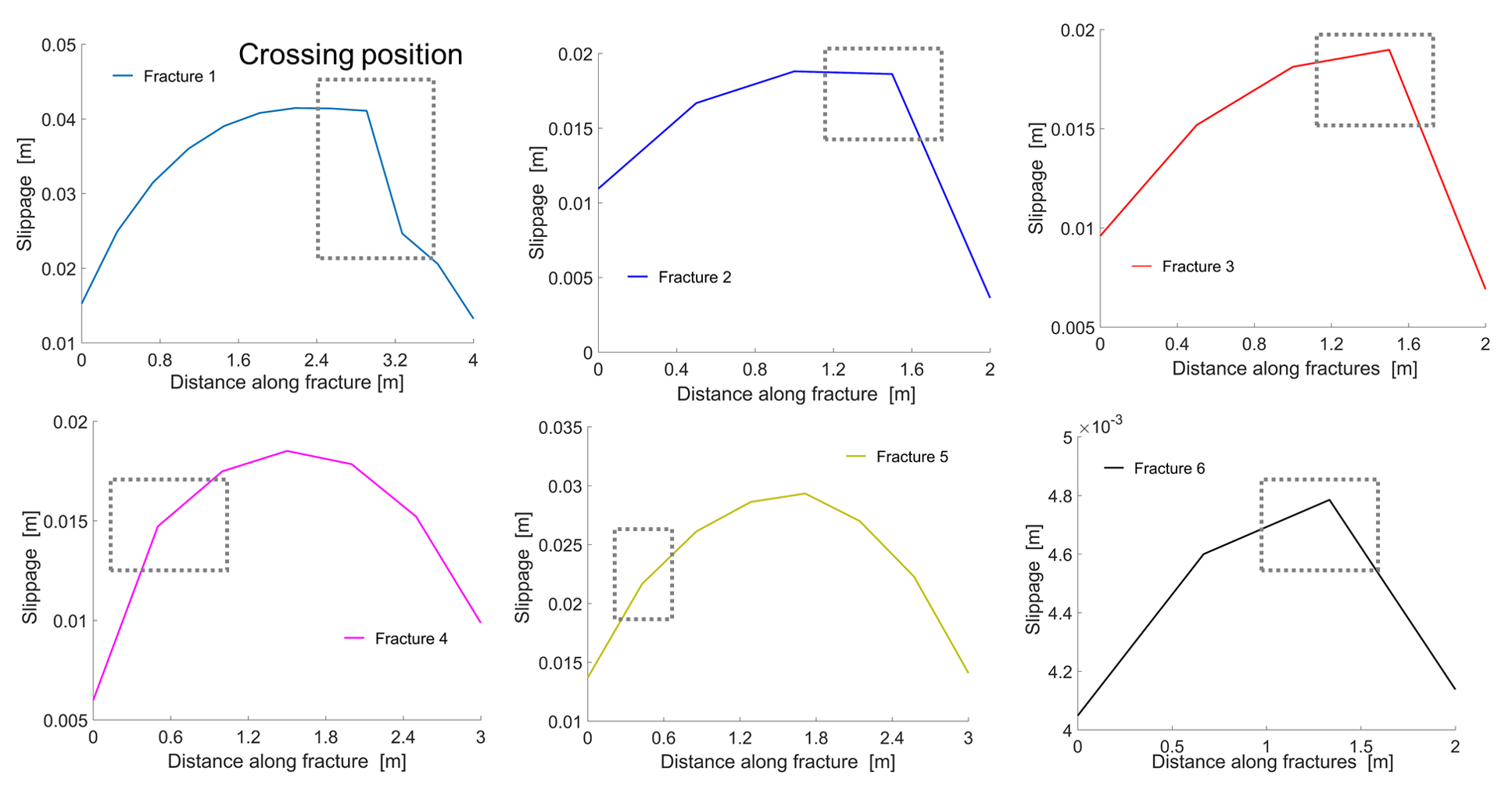}
\caption{ The slippage on contact surface of the multiple crossing fractures model. The top row for Fracture $1 \sim 3$. The bottom row for Fracture $4 \sim 6$. Note that the portions marked by dash boxes imply the crossing positions. }
\label{fig:fig_Test6_curves}
\end{figure}

\begin{figure}[H]
\centering
\includegraphics[width=14.5cm]{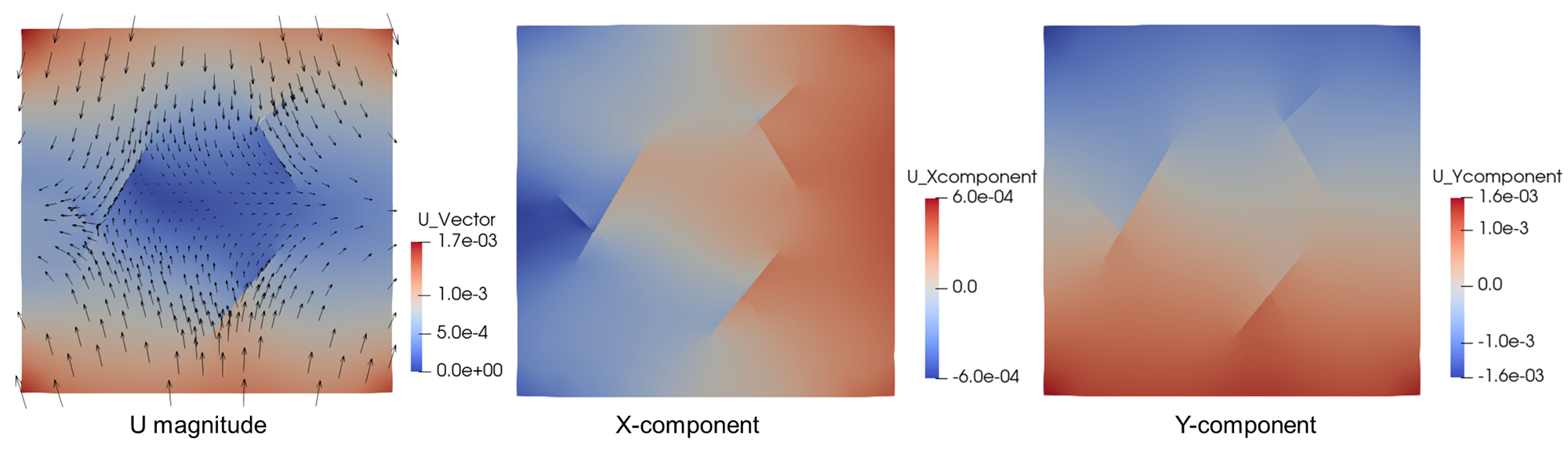}
\caption{ Deformation of multiple crossing fractures model. The vector field (left) and the displacement components (middle and right). }
\label{fig:fig_multiCrossingFrac_countor}
\end{figure}

The schematic of multiple crossing fractures is displayed in Fig. \ref{fig:fig_Test5and6_GeometricalModel}b. A set of multiple crossing fractures is simulated in this test. The simulation results show the stability of the proposed numerical method when it treats several crossing fractures. The displacement field is shown in Fig. \ref{fig:fig_multiCrossingFrac_countor}. It appears that the displacement vector arrows have some shifts around the fractures, as seen in the left picture. Fig. \ref{fig:fig_Test6_curves} shows the slippage on each fracture of this fractured domain, which are labeled in Fig. \ref{fig:fig_Test5and6_GeometricalModel}b. The shape of each curve is similar with the curve in the test of single fracture, i.e. the parabolic mode, as discussed in Section \ref{subsection:Benchmark1}. The key feature among the six fractures is the jump displacement at the crossing positions produced by the multi-intersected fractures.

\begin{figure}[H]
\centering
\includegraphics[width=14.5cm]{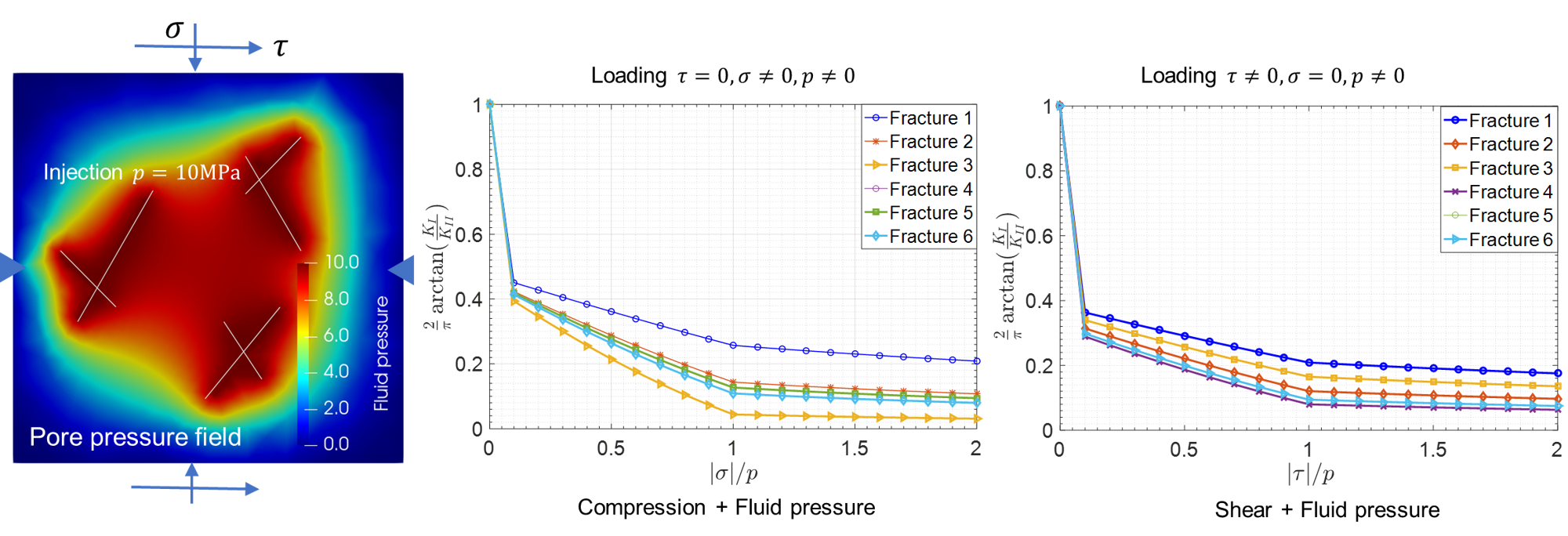}
\caption{ Variation of the normalized stress intensity factor (SIF) ratio under different loading ratios (middle and right). An example of pore pressure field when $p=10$ MPa (left). Fracture numbers are labeled in Fig. \ref{fig:fig_Test5and6_GeometricalModel}. }
\label{fig:fig_SIFCurves}
\end{figure}

Fig. \ref{fig:fig_SIFCurves} displays the relation of the normalized stress intensity factor (SIF) ratio versus the loading ratio. Concerning the pure mode I or mode II, the SIFs are $K_I$ and $K_{II}$, respectively. 
The normalized SIF ratio is commonly used to evaluate the domination of mode-I and mode-II failure \cite{Wang2019a}. It is defined as $({2}/{\pi})\arctan(K_I / K_{II})$. We test two different mixed mode loadings, which are the compression or tension combines the fluid pressure imposed on internal surface of fractures. It appears that the fluid pressure induces a larger opening of fractures, also know as the mode I failure. In contrast, the compression and tension loadings lead to a shear failure (mode II failure). With the increase of $\sigma$ or $\tau$, the normalized SIF ratio shows a decrease tendency. The shear failure becomes the domination compared to fracture opening (mode I), especially when the loading ratios $|\sigma / p|$ and $\tau /p$ are greater than 1.

\section{Conclusions}
\label{section:Conclu}

In this work, a improved mixed-FE scheme is proposed using the Lagrange multiplier method in the framework of constrained variational principle, which has the capability to handle frictional contact and shear failure of the multiple crossing fractures under mixed mode loadings.

The contact model of multiple fractures is generalized based on a rigorous mathematical formulation of contact theory. Both the isolated fractures and crossing fractures can be treated in the framework. The compression would lead to the activation of contact constraints, while the tensile loading remains the standard treatment of elastic deformation. The additional constraints of fracture contact are integrated into the constrained variational principle (CVP) through Lagrange multipliers. In this way, an unified formulation is derived, in which the slippage, opening and contact traction on the fracture surface can be calculated by the resulted saddle-point algebraic system.

Galerkin finite element approximation is used to discretize the governing equations. The discretized contact pairs along each of the fractures are defined to capture the interaction between the two sides of a fracture. To handle multiple crossing fractures, the crossing contact pairs are devised to overcome the unphysical scenario at the crossing position. 
The two primary unknowns of the system, namely displacement and Lagrange multiplier, are solved by an iteration method. For a robust numerical scheme, the precondition technique is introduced to handle the saddle-point algebraic system, which leads to a preconditioned mixed-FE scheme. 
The contact system is then resolved using a monolithic iteration strategy, motivated by the high nonlinear property of the governing equations.

A series of numerical tests is conducted to study the contact behaviors of single- and multiple crossing fractures based on the proposed mixed-FE scheme. First, the benchmark study is presented to verify the numerical results. Later, the contact on fractures is captured by mixed mode loading. Two tests with complex geometry are studied. Especially, the effects of crossing fractures on the deformation field can be observed in the calculated results, in which the slippage, opening and stress intensity factor are analyzed with different conditions.

\section*{Declaration of Competing Interest}

The authors declare that they have no known competing financial
interests or personal relationships that could have appeared to influence
the work reported in this paper.

\section*{Acknowledgements}

Financial support of NWO-TTW ViDi grant (project ADMIRE 17509) is acknowledged. 
Luyu Wang gratefully acknowledges the financial support by China Scholarship Council (No. 201904910310). 
The authors would like to thank all the members of the DARSim (Delft Advanced Reservoir Simulation) research group and ADMIRE project members at TU Delft, for the fruitful discussions.



\appendix
\section{Weak forms derived by the constrained variational principle}
\label{appendix_weakForms}

Following the formulation in Section \ref{section:SecDiscretization}, we have the weak form if the fractures are stick: 
\begin{footnotesize}
\begin{equation}
\begin{aligned}
\delta \Pi^{CL}_{u} \left( \bm{u}^f ,\bm{\lambda} \right) 
=\sum_{i=1}^{N^f} 
\int_{\Gamma_{in,i}^{stick}}  \bm{\lambda}^T \delta \bm{C} \left( \bm{u}^{f} , \bm{t}^{f} \right) d\Gamma 
= 
\sum_{i=1}^{N^f} \int_{\Gamma_{in,i}^{stick}}  \delta [\![ \hat{\bm{u}}^f ]\!]^T  \bm{\lambda}  d\Gamma  
\end{aligned}
\label{eq:WeakForm_State_Disp}
\end{equation}
\end{footnotesize}
with the relative displacement vector $[\![ \hat{\bm{u}}^f ]\!] = [[\![ {{u}}_N^f ]\!]  \;\; [\![ {{u}}_T^f ]\!]]^T$ in the local coordinate system attached to a certain fracture surface $\Gamma_{in,i}$, as shown in Fig. \ref{fig:fig_LocalSystem}. 
Meanwhile, Eq. (\ref{eq:VariationalPrinciple_Lambda}) is written as: 
\begin{footnotesize}
\begin{equation}
\begin{aligned}
\delta \Pi^{CL}_{\lambda} \left( \bm{u}^f ,\bm{\lambda} \right)
&=
\sum_{i=1}^{N^f} 
\int_{\Gamma_{in,i}^{stick}} \delta \bm{\lambda}^T \bm{C} \left( \bm{u}^{f} , \bm{t}^{f} \right) d\Gamma 
=
\sum_{i=1}^{N^f} 
\int_{\Gamma_{in,i}^{stick}} \delta \bm{\lambda}^T  \left( [\![ \bm{\hat{u}}^f ]\!] + \bm{g} \right) d\Gamma 
\end{aligned}
\label{eq:WeakForm_State_Lambda}
\end{equation}
\end{footnotesize}
with the initial gap $\bm{g} = [g_N \;\; 0]^T$ between $+$ side and $-$ side, as indicated in Eq. (\ref{eq:KKTCondition_Normal}). Normally, $g_N$ equals either 0 or a small value such as $10^{-3} \sim 10^{-4} {\rm m}$.

Following the minimization of energy functional \cite{Lanczos2012,Reddy2017}, the constrained variational principle with Lagrange multiplier requires: 
\begin{footnotesize}
\begin{equation}
\begin{aligned}
\int_{\Omega} \delta \bm{\varepsilon}^T \bm{\sigma} d\Omega
- \int_{\Omega} \delta \bm{u}^T \bm{f} d\Omega 
- \int_{{{\Gamma}}_{ex}^N} \delta \bm{u}^T \bm \bar{\bm{t}}^{ex} d\Gamma 
+ 
\sum_{i=1}^{N^f}  \int_{\Gamma_{in,i}^{stick}}  \delta [\![ \hat{\bm{u}}^f ]\!]^T  \bm{\lambda}  d\Gamma 
= 0
\\
\sum_{i=1}^{N^f} 
\int_{\Gamma_{in,i}^{stick}} \delta \bm{\lambda}^T   [\![ \bm{\hat{u}}^f ]\!]   d\Gamma 
= -  \sum_{i=1}^{N^f} \int_{\Gamma_{in,i}^{stick}}  \delta \bm{\lambda}^T  \bm{g}  d\Gamma 
\end{aligned}
\label{eq:WeakFormState}
\end{equation}
\end{footnotesize}

The two equations in system Eq. (\ref{eq:WeakFormState}) are independent to each other. We calculate the unknowns, displacement $\bm{u}$ and Lagrange multiplier $\bm{\lambda}$, in a framework of numerical scheme, namely the mixed-FE scheme.

For slip sate, the first-order variation of energy functional with respect to displacement is given by: 
\begin{small}
\begin{equation}
\begin{aligned}
\int_{\Omega} \delta \bm{\varepsilon}^T \bm{\sigma} d\Omega
&- \int_{\Omega} \delta \bm{u}^T \bm{f} d\Omega 
- \int_{{{\Gamma}}_{ex}^N} \delta \bm{u}^T \bm \bar{\bm{t}}^{ex} d\Gamma 
+ 
\sum_{i=1}^{N^f} \left( \int_{\Gamma_{in,i}^{slip}} [\![ \delta u_N^f ]\!] \lambda_N  d\Gamma 
\right. \\& \left.
+ \int_{\Gamma_{in,i}^{slip}} [\![ \delta u_T^f ]\!] c \; sign  d\Gamma 
- \int_{\Gamma_{in,i}^{slip}} [\![ \delta u_T^f ]\!] \lambda_N sign \tan{\varphi} d\Gamma 
\right)
= 0 
\end{aligned}
\label{eq:WeakFormSlip_Disp}
\end{equation}
\end{small}
and the variation with respect to Lagrange multiplier is given by: 
\begin{small}
\begin{equation}
\begin{aligned}
\sum_{i=1}^{N^f} 
\int_{\Gamma_{in,i}^{slip}} \delta \lambda_N \left(    [\![ {{u}}_N^f ]\!]   -  [\![ {{u}}_T^f ]\!] sign\tan \varphi   \right) d\Gamma 
= -  \sum_{i=1}^{N^f} \int_{\Gamma_{in,i}^{slip}}  \delta {\lambda}_N  g_N  d\Gamma 
\end{aligned}
\label{eq:WeakFormSlip_Lambda}
\end{equation}
\end{small}
which holds in the scenario of multi-fractures indicated by the number of fractures $N^f$. Other notations are defined in the preceding sections.

\section{Discretized forms of the mixed-FEM scheme}
\label{appendix_discretizedForms}

Following the formulation in Section \ref{section:SecSolutionStrategy}, the components in Eq. (\ref{eq:SaddlePointSystem}) are expressed as: 
\begin{small}
\begin{equation}
\begin{aligned}
& \bm{R}^u = \bm{K}^{uu} \bm{U} + \bm{C}^{\lambda u} \bm{\Lambda} - \left( \bm{F} - \bm{F}^{slip} \right) \\
& \bm{R}^{\lambda} = \left[\bm{C}^{\lambda u}\right]^T \bm{U} + \bm{g} \\
& \bm{C}^{\lambda u} = \bm{K}^{\lambda u} + \left( \bm{K}^{\lambda u,N} - \bm{K}^{\lambda u, T} \right)\\
\end{aligned}
\label{eq:Blocks1}
\end{equation}
\end{small}
and the sub-blocks: 
\begin{small}
\begin{equation}
\begin{aligned}
& \bm{K}^{uu} = \int_{\Omega} \bm{B}^T \bm{D} \bm{B} d\Omega  \\
& \bm{K}^{\lambda u} = \int_{\Gamma_{in}^{stick} \cup \Gamma_{in}^{slip}} \bm{G}^T \bm{S} \bm{N}^{\lambda}  d\Gamma \\
& \bm{K}^{\lambda u,N} =  \int_{ \Gamma_{in}^{slip}} \bm{G}^T \bm{n}^f \bm{N}^{\lambda}_N  d\Gamma  \\
& \bm{K}^{\lambda u,T} = \int_{ \Gamma_{in}^{slip}} \bm{G}^T \bm{m}^f \bm{N}^{\lambda}_N  sign \tan{\varphi} d\Gamma \\
& \bm{F} = \int_{\Omega} \left[ \bm{N}^u \right]^T \bm{f} d\Omega + \int_{\Gamma_{ex}^N}  \left[ \bm{N}^u \right]^T \bar{\bm{t}}^{ex}d\Gamma \\
& \bm{F}^{slip} = \int_{\Gamma_{in}^{slip}} \bm{G}^T \bm{m}^f c \; sign d\Gamma
\end{aligned}
\label{eq:Blocks2}
\end{equation}
\end{small}
where the blocks $\bm{F}^{slip}$, $\bm{K}^{\lambda u,N}$ and $\bm{K}^{\lambda u,T}$ are related to the contributions of slip contact. $\bm{K}^{\lambda u}$ is related to the co-effects by stick and slip states. $sign$ is defined in Eq. (\ref{eq:IndicatorFunction}).

The system of Eqs. (\ref{eq:SaddlePointSystem}), (\ref{eq:Blocks1}) and (\ref{eq:Blocks2}) will be constructed at each iteration step. The iteration process is controlled by the termination criterion: 
\begin{equation}
\left\|  \mathbf{R}^{n+1}  \right\|_2  \;  <  \; \epsilon 
\end{equation}
with a user-defined threshold $\epsilon $. Normally, it is set to a default value $10^{-4}$. All notations are defined in the preceding sections.

Eqs. (\ref{eq:SaddlePointSystem}), (\ref{eq:Blocks1}) and (\ref{eq:Blocks2}) render one of the most innovations and a core contribution of the presented work, which include the unified consideration of stick/slip/open states for the contact response of crossing fractures.

\section{Expanded forms of the preconditioner}
\label{appendix_preconditioner}

Following the formulation in Section \ref{subsection:SubsecPreconditioning}, the expanded forms are written as: 
\begin{equation}
\bm{A}^2 = 
\begin{bmatrix}
\sum_{j=1}^{n_{c}^E} (E_{1j})^2 & \cdots & \cdots	\\
\cdots & \sum_{j=1}^{n_{c}^E} (E_{2j})^2 & \cdots \\
\vdots  &  \ddots &  \vdots  \\ 
\cdots & \cdots & \sum_{j=1}^{n_{c}^E} (E_{n_{r}^Ej})^2  
\end{bmatrix}
\end{equation}
\begin{equation}
\bm{B}^2 = 
\begin{bmatrix}
\sum_{j=1}^{n_{c}^D} (D_{1j})^2 & \cdots & \cdots	\\
\cdots & \sum_{j=1}^{n_{c}^D} (D_{2j})^2 & \cdots \\
\vdots  &  \ddots &  \vdots  \\ 
\cdots & \cdots & \sum_{j=1}^{n_{c}^D} (D_{n_{r}^Dj})^2  
\end{bmatrix}
\end{equation}
with the numbers of columns $n_{c}^D$ and rows $n_{r}^D$ for $\bm{D}$, $n_{c}^E$ and $n_{r}^E$ for $\bm{E}$. Other notations are defined in the preceding sections. 
This operation provides a way to improve the numerical quality of the Jacobian. It is efficient and easy to implement.

\section{The monolithic-updated contact algorithm}
\label{appendix_algorithm}

\begin{algorithm}[htb]
\footnotesize
\caption{ The monolithic-updated contact algorithm on unstructured grids}

\label{Algorithm1}
\begin{algorithmic}[1]

\FOR{each $i \in [1,N^f]$}
\STATE Fracture $\Gamma_{in,i}$ is determined by its coordinate $\left( x_i, y_i \right)$
\STATE Obtain the node connectivity and contact pairs on unstructured grid
\ENDFOR

Enter the procedure of contact state update

\FOR{each $t \in [0,t_{max}]$}
\STATE At time step $t$, calculate the current contact state of contact pairs
\label{code:Algorithm1:a1}
\STATE Enter the iteration process
    \FOR{each $\nu \in [0,\nu_{max}]$}
    \STATE At iteration $\nu$, construct the system Eq. (\ref{eq:SaddlePointSystem}) based on current contact state
\label{code:Algorithm1:a2}
    \STATE Obtain the results of contact state at iteration $\nu+1$, as indicated in Eq. (\ref{eq:SystemMatrixForm})
    \STATE Update the unknowns $\bm{U}^{\nu+1}=\bm{U}^{\nu}+\delta \bm{U}^{\nu+1}$ and $\bm{\Lambda}^{\nu+1}=\bm{\Lambda}^{\nu}+\delta \bm{\Lambda}^{\nu+1}$
    \ENDFOR
\STATE If the contact state is satisfied with contact constraints then go to the updated time step $t=t+1$, and refer to Line \ref{code:Algorithm1:a1}
\STATE Otherwise, go to Line \ref{code:Algorithm1:a2} until the calculated contact state is satisfied with contact constraint

\IF {$t=t_{max}$ }
\STATE Terminate the program
\ENDIF

\ENDFOR

\end{algorithmic}
\end{algorithm}




\end{document}